\documentclass[aop,MSNbibl,citesort,dvips]{arximspdf}

\doi{10.1214/09-AOP475}
\volume{38}
\issue{1}
\pubyear{2010}
\firstpage{207}
\lastpage{233}

\makeatletter

\newtheorem{theorem}{Theorem}
\newtheorem{coro}[theorem]{Corollary}
\newproclaim{remark}[theorem]{Remark}
\newtheorem{lemma}[theorem]{Lemma}
\newtheorem{proposition}[theorem]{Proposition}

\makeatother

\begin{document}
\begin{frontmatter}

\title{The $\Lambda$-coalescent speed of coming down from infinity}
\runtitle{The speed of $\Lambda$-coalescents}

\begin{aug}
\author[A]{\fnms{Julien} \snm{Berestycki}},
\author[B]{\fnms{Nathana\"{e}l} \snm{Berestycki}\corref{}\ead[label=e1]{N.Berestycki@statslab.cam.ac.uk}} and
\author[C]{\fnms{Vlada} \snm{Limic}\thanksref{t1}}
\runauthor{J. Berestycki, N. Berestycki and V. Limic}
\affiliation{Universit\'{e} Pierre et Marie Curie-ParisVI, Cambridge University and C.N.R.S.}
\address[A]{J. Berestycki\\
Laboratoire de Probabilit\'{e}s\\
\quad et Mod\`{e}les Al\'{e}atoires / UMR 7599\\
Universit\'{e} Pierre et Marie Curie-Paris VI\\
4 Place Jussieu\\
F-75252 Paris Cedex 05\\
France} 
\address[B]{N. Berestycki\\
Statistical Laboratory\\
 University of Cambridge\\
Wilberforce Rd.\\
 Cambridge CB3 0WB\\ United Kingdom\\
 \printead{e1}}
\address[C]{V. Limic\\
 Universit\'{e} de Provence\\
Technop\^{o}le de Ch\^{a}teau-Gombert\\
UMR 6632, LATP, CMI\\
39, rue F. Joliot Curie\\
13453 Marseille, cedex 13\\ France}
\end{aug}

\thankstext{t1}{Supported in part by
NSERC Discovery grant and by Alfred P. Sloan Research Fellowship.}

\pdfauthor{Julien Berestycki, Nathanael Berestycki, Vlada Limic}

\received{\smonth{7} \syear{2008}}
\revised{\smonth{2} \syear{2009}}

%
\begin{abstract}
Consider a $\Lambda$-coalescent that
comes down from infinity (meaning that
it starts from a configuration containing
infinitely many blocks at time $0$, yet
it has a finite number $N_t$ of blocks at any positive time
$t>0$). We
exhibit a deterministic function $v\dvtx(0,\infty)\to(0,\infty)$
such that $N_t/v(t)\to1$, almost surely, and in $L^p$ for any $p\geq1$,
as $t\to0$.
Our approach relies on a novel martingale technique.
\end{abstract}

%

\end{frontmatter}

\section{Introduction}\label{sec1}

Various natural population
genetics models lead to a representation of the genealogical tree
by a process called Kingman's coalescent \cite{king82,king82b}.
Kingman's coalescent is a Markov process which can be informally
described as follows:
in a fixed sample of $n$ individuals from the
population, each pair of ancestral lineages coalesces at rate 1.

In population genetics, one uses the above process to quantify
polymorphism in a homogeneously mixing population under neutral
evolution.
However, there is some evidence that for
modeling evolution of marine populations (see, e.g., \cite{marine}),
the use of coalescent processes which allow \textit{multiple
collisions} is more appropriate than that of Kingman's coalescent
where only pairs of blocks can merge at any given time.
Similarly,
multiple collisions are natural for modeling evolution of viral populations,
where natural selection plays a very strong role.
They also emerge in the fine-scale mapping of disease loci \cite{morrisetal}.

A suitable family of mathematical models has been introduced and
studied by Pitman \cite{pit99} and Sagitov \cite{sag99}
under the name \textit{$\Lambda$-coalescents} or
\textit{coalescents with multiple collisions}.
We postpone
the precise definitions of these processes until the next section.

Let $N^\Lambda\equiv N:=(N_t, t \geq0)$ be the number of blocks process
corresponding to a particular $\Lambda$-coalescent process.
In view of applications,
we concentrate on $\Lambda$-coalescents such that $P(N_t<\infty,
t>0)=1$ and
$\lim_{t\to0+} N_t=\infty$ (here $N$ is really an entrance law).
This property is typically referred to as coming down from infinity
(see Section \ref{S:Lacoal}
for a formal definition).
It is important to understand the nature of divergence of $N_t$
as $t$ decreases to $0$.
In the current paper, our goal is
to exhibit a function $v\dvtx(0,\infty)\to(0,\infty)$ such that
\[
\lim_{t\to0}
\frac{N_t}{v(t)}= 1 \qquad\mbox{almost surely.}
\]
We call any such $v$ the \textit{speed} of coming down from infinity
(speed of CDI)
for the corresponding $\Lambda$-coalescent.
Note that the limit above is
in fact the limit as $t\to0+$; from now on we always write $t \to0$.
The exact form of the function $v$ is implicit and somewhat
technical (see Theorem \ref{Tsmalltime} for the precise statement).
However, in many situations of interest, one can find a simpler
function $g(t)$, often a power in $t$, such that $g(t)/v(t)\to1$,
and therefore $N_t/g(t)\to1$, as $t\to0$. Then we also refer to
$g$ as the speed of CDI for the corresponding coalescent. As
mentioned above, Kingman's coalescent is the simplest
$\Lambda$-coalescent. In particular, one can quickly find its speed of
CDI by considering the ``time-reversed'' process. Analogous
time-reversals for general $\Lambda$-coalescents seem to be difficult to
grasp. The speed of CDI was recently determined for Beta-coalescents
and their ``perturbations'' in Berestycki, Berestycki and Schweinsberg
\cite{bbs2} and \cite{bbs1}
and Bertoin and Le Gall \cite{blg3} (where convergence
is established in probability). See also the comment following the
statement of Theorem \ref{Tsmalltime} below.

With the above biology
motivation in mind, there is a strong interest in
understanding
(see, e.g., \cite{dong,eldwak,mohle})
analogues of Ewens' sampling formula for $\Lambda$-coa\-lescents.
It seems that only Kingman's coalescent
allows for an exact solution
(see, e.g., \cite{durrettDNA} or \cite{ewens})
while in the general case, one should aim for good approximations.
The only previous detailed analysis of this kind was carried out
in \cite{bbs1} and \cite{bbs2}
for the special case of Beta-coalescents.
The above result can be viewed as the first step towards analogous
understanding of the general $\Lambda$-coalescent case.

In a parallel work \cite{bbl2} we discuss
the consequence of our main results to the problem of
quantifying polymorphism in a population whose genealogy is driven by a
coalescent
with multiple collisions.
In the same paper, we will describe a general connection between the
small-time asymptotics of $\Lambda$-coalescents
and continuous random trees and their associated continuous-state
branching processes as well as
generalized Fleming--Viot processes.
These connections enable one to guess the form of function $v(t)$,
and they imply the convergence in probability of
the quantity $N_t/v(t)$ which is of interest under certain technical conditions.
They can also be useful in determining the power law
order of growth of $v$ as $t\to0$.

To the best of our knowledge, the martingale analysis in the current
context is novel.
We believe that it is of independent interest.
Although similar in spirit,
our setting is different from the general setting of Darling and Norris
\cite{darnor}.
For their technique to apply,
it is necessary to start with good bounds on the
accumulated absolute difference of the ``drifts'' of the Markov chain
and the
solution to the corresponding differential ``fluid-limit'' equation.
Here it seems difficult to obtain such bounds. However, it is possible
to work directly
[cf. the local martingale $M_z'$ from (\ref{Elocmart})] 
with the accumulated (nonabsolute) difference of the drifts
in order to obtain sufficiently good asymptotic estimates.

The rest of the paper is organized as follows. Section \ref{S:prelim}
contains definitions and notations.
The main results are stated in Section \ref{S:results} and are proved
in Section \ref{S:martingale}, with
some technical estimates postponed until the \hyperref[ap:binom]{Appendix}.

\section{Definitions and preliminaries}
\label{S:prelim}
\subsection{Notation}
\label{S:notation} We recall some standard notation, and introduce
additional notation to simplify the exposition.
\begin{itemize}
\item[{$ $}]
Denote the set of real (resp. rational) numbers by $\mathbb{R}$
(resp. $\mathbb{Q}$) and set $\mathbb{R}_+ = (0,\infty)$.
For $a,b\in\mathbb{R}$, denote by $a \wedge b$ (resp. $a \vee b$) the minimum
(resp. maximum)
of the two numbers.
\item[{$ $}]
Let $\mathbb{N}:=\{1,2,\ldots\}$ and let $\mathcal{P}$ be the set of
partitions of $\mathbb{N}$.
Furthermore, for $n\in\mathbb{N}$ denote by
$\mathcal{P}_n$ the set of partitions of $[n]:=\{1,\ldots, n\}$.
\item[{$ $}]
If $f$ is a function, defined in a left-neighborhood $(s-\varepsilon
,s)$ of a
point $s$,
we denote by $f(s-)$ the left limit of $f$ at $s$.
\item[{$ $}]
Given two functions $f,g\dvtx\mathbb{R}_+\to\mathbb{R}_+$, write
$f=O(g)$ if $\limsup f(x)/g(x) <\infty$,
$f=o(g)$ if $\limsup f(x)/g(x) =0$,
and $f\sim g$ if $\lim f(x)/g(x) =1$.
The point at which the limits are taken
might vary, depending on the context.
\item[{$ $}]
If $X$ and $Y$ are two random objects, we write $X\stackrel{d}=Y$ to
indicate their
equivalence in distribution.
As usual, convergence in distribution will be denoted by
$\Rightarrow$ symbol.
\item[{$ $}]
If $\mathcal{F}=(\mathcal{F}_t,t\ge0)$ is a filtration, and $T$ is a
stopping time
relative to $\mathcal{F}$, denote by $\mathcal{F}_T$ the standard filtration
generated by $T$ (see, e.g., \cite{durrett}, page 389).
\item[{$ $}]
For $\nu$ a finite or $\sigma$-finite measure, denote the
support of $\nu$ by $\operatorname{supp} (\nu)$.
\end{itemize}

\subsection{$\Lambda$-coalescents}
\label{S:Lacoal}
Let $\Lambda$ be a finite measure on $[0,1]$.
The $\Lambda$-coalescent is a Markov process $(\Pi_t,t\ge0)$ with values
in $\mathcal{P}$ (the set of partitions of $\mathbb{N}$),
characterized as follows. If $n \in\mathbb{N}$,
then the restriction
$(\Pi^{(n)}_t,t\ge0)$
of $(\Pi_t,t\ge0)$ to
$[n]$ is a Markov chain, taking values in $\mathcal{P}_n$,
with a following dynamics:
whenever $\Pi^{(n)}_t$ is a partition consisting of $b$ blocks, the
rate at which a given $k$-tuple of its blocks merges is
%
%
\begin{equation}
\label{rate_coal}
\lambda_{b,k}=\int_{[0,1]}x^{k-2}(1-x)^{b-k}\Lambda(dx).
\end{equation}
Note that mergers of several blocks into one block are possible, but
multiple mergers do not occur simultaneously. For a generalization
of $\Lambda$-coalescents where multiple mergers are possible, see
Schweinsberg \cite{sch2}. For a generalization of
$\Lambda$-coa\-lescents to spatial (not a mean-field) setting, see
Limic and Sturm \cite{ls}.

We will quote here several basic
properties of the $\Lambda$-coalescent, and refer the reader to
Pitman \cite{pit99} for details and additional analysis.
When $\Lambda(\{0\})=0$, the corresponding
$\Lambda$-coalescent can be constructed via a Poisson point process in
the following way. Let
%
%
\begin{equation}
\label{DPPPpi}
\pi(\cdot) = \sum_{i \in\mathbb{N}} \delta_{t_i,x_i}(\cdot)
\end{equation}
be a Poisson point process on $\mathbb{R}_+ \times(0,1)$ with the
intensity measure
$dt \otimes\nu(dx)$ where $\nu(dx)= x^{-2}\Lambda(dx)$.
Each atom $(t,x)$ of $\pi$ influences the evolution of the process
$\Pi$ as
follows:
for each block of $\Pi(t-)$, flip a coin with probability of
heads equal to $x$; all the blocks corresponding to
coins that come up ``head'' then merge immediately
into one single block while all other blocks remain unchanged.
Note that in order to make this construction rigorous,
one first considers the restrictions
$(\Pi^{(n)}(t),t\ge0)$, since the measure $\nu(dx)=x^{-2}\Lambda
(dx)$ may have infinite total mass.

We next recall a remarkable
property
of $\Lambda$-coalescents.
Let $E$ be the event that for all $t>0$ there are infinitely many
blocks, and let $F$ be the event that for all $t>0$ there are only
finitely many blocks.
Pitman \cite{pit99} showed
that, if $\Lambda(\{1\})=0$,
only the following two types of behavior are possible, depending on the
measure~$\Lambda$:
either $P(E)=1$
or $P(F)=1$.
When $P(F) = 1$, the process $\Pi$ is said to \textit{come down
from infinity}. For instance, Kingman's coalescent comes down from
infinity, while if $\Lambda(dx)=dx$ is the uniform measure on
$(0,1)$, then the corresponding $\Lambda$-coalescent does not come
down from infinity.
This particular $\Lambda$-coalescent was discovered by Bolthausen and
Sznitman \cite{bosz98}
in connection with spin glasses.

A necessary and sufficient condition for a $\Lambda$-coalescent to
come down from infinity was given by Schweinsberg
\cite{sch1}: define
\[
\gamma_b = \sum_{k=2}^b (k-1)\pmatrix{
b \cr k}
\lambda_{b,k},
\]
then the $\Lambda$-coalescent comes down from infinity if and only if
$\sum_{b=2}^{\infty} \gamma_b^{-1} <\infty$.

Recently, Bertoin and Le Gall \cite{blg3} observed that
this condition is equivalent to the following requirement: define
%
%
\begin{equation}
\label{D:psi}
\psi_{\Lambda}(q)\equiv\psi(q):= \int_{[0,1]}(e^{-qx}-1+qx)\nu(dx),
\end{equation}
where $\nu(dx)=x^{-2}\Lambda(dx)$,
then
%
%
\begin{equation}\label{EequivaleG}
\sum_{b=2}^{\infty} \gamma_b^{-1} <\infty
\quad\mbox{if and only if}\quad \int_a^{\infty} \frac{dq}{\psi(q)}
<\infty,
\end{equation}
where the right-hand side is finite for some (and then automatically
for all) $a >0$.
Somewhat remarkably, the divergence rate function $v$
is given [cf. definition (\ref{Ev}) in the next section] in terms of
the right-hand side in
(\ref{EequivaleG}).
The condition (\ref{EequivaleG}) is well known
in the L\'{e}vy processes literature as the Grey's criterion for
extinction of the underlying continuous-state branching process. We
refer the reader to \cite{bbl2} for further explanation of the above
connections.

\section{Main results}
\label{S:results}
Let $\Lambda$ be a finite measure
on $[0,1]$, and let $(\Pi_t,t\ge0)$ be a $\Lambda$-coalescent.
Without loss of generality, we may, and will, henceforth assume that
$\Lambda$ is a
probability measure, that is,
%
%
\begin{equation}
\label{unit_mass}
\Lambda[0,1]=1.
\end{equation}
Indeed, a scaling of the total mass of $\Lambda$ by a
constant factor will induce the scaling of the speed of evolution (and
therefore, that of
coming down from infinity) by the same factor, and
the speed of CDI $v$ from (\ref{Ev}) below will
scale in the same way.

To each such measure $\Lambda$ we associate a function $\psi$ defined in
(\ref{D:psi}).
Moreover, for a probability measure $\tilde\Lambda$
of the form $\tilde\Lambda=(1-c)\Lambda+c \delta_0$, where $\Lambda
$ has
no atom at $0$,
we may rewrite
as
%
%
\begin{equation}
\label{Epsi_atom} \psi_{\tilde\Lambda}(q) = \frac{c}{2} q^2 +
(1-c)\int_{[0,1]}(e^{-qx}-1+qx)\nu(dx).
\end{equation}
Note that if $c=1$
we retrieve the Kingman coalescent, whose small-time behavior is well
understood.
Henceforth we assume that $c<1$.

When $\psi$ is such that
the integral in (\ref{EequivaleG}) is finite, or equivalently,
when the corresponding $\Lambda$-coalescent comes down from infinity, we
can define
%
%
\begin{equation}\label{Eu}
u_{\psi}(t)\equiv u(t):=\int_t^{\infty} \frac{dq}{\psi(q)} \in
(0,\infty),\qquad t>0,
\end{equation}
and its c\`{a}dl\`{a}g inverse
%
%
\begin{equation}
\label{Ev} v_\psi(t)\equiv v(t) :=\inf\biggl\{s>0\dvtx\int_s^{\infty}
\frac1{\psi(q)}\,dq
<t \biggr\}, \qquad t>0.
\end{equation}
Denote by
$(N^\Lambda(t), t\geq0)=(N^\Lambda_t,t\geq0)$ the number of blocks process
for the
$\Lambda$-coalescent $(\Pi(t),t\geq0)$.
The first main result of this paper is following theorem.
\begin{theorem} \label{Tsmalltime}
%
%
\begin{equation}\label{nblocks}
\lim_{t\to0}\frac{N^\Lambda(t)}{v_\psi(t)} = 1\qquad\mbox{almost surely}.
\end{equation}
\end{theorem}

Note that if $\Pi$ does not come from infinity, both $N^\Lambda
_t=N^\Lambda
(t)=\infty$, for all $t\geq0$,
almost surely, and
the formal definition (\ref{Ev}) yields $v_\psi\equiv\infty$,
so (\ref{nblocks}) extends trivially if $\infty/\infty=1$.

We next comment on some special cases of Theorem \ref{Tsmalltime}. When
$\Lambda=\delta_{0}$, we have $v(t)= 2/t$, and we recover the
well-known result that for Kingman's coalescent, the number of
blocks is almost surely asymptotic to $2/t$. Another interesting
case occurs when $\Lambda$ has the $\operatorname{Beta}(2-\alpha,\alpha)$
distribution for some $1<\alpha<2$. That is,
%
%
\begin{equation}\label{D:beta}
\Lambda(dx) = \frac{1}{\Gamma(2-\alpha) \Gamma(\alpha)}
x^{1-\alpha} (1-x)^{\alpha- 1} \,dx.
\end{equation}
Here it is not hard to see that
$\psi(q)\sim c_1 q^{\alpha}$ as $q\to\infty$, and thus
that
\[
v(t) \sim c_2t^{-1/(\alpha-1)}\qquad\mbox{as } t\to0,
\]
where $c_1=(\Gamma(\alpha) \alpha(\alpha-1))^{-1}$ and
$c_2=(\alpha\Gamma(\alpha))^{-1/(\alpha-1)}$. In fact these
calculations can easily be generalized to the case where $\Lambda$
is regularly varying near 0 with index $1<\alpha<2$. In this
case, Theorem \ref{Tsmalltime} strengthens
Lemma 3 in \cite{blg3}.

However, we emphasize that the most delicate case of the above
theorem occurs when the
measure $\Lambda$ ``wildly oscillates'' in any neighborhood of 0.
An example of such a measure is constructed in
the appendix of \cite{bbl2}.
It illustrates potential difficulties in the analysis of
functions $\psi$, $u$ or $v$ directly.

With a bit more work, we obtain as the second main result
an analogue to Theorem \ref{Tsmalltime} in terms
of convergence of moments.
\begin{theorem} \label{Tsmalltime-mom}
For any $d \in[1, \infty)$,
%
%
\begin{equation}
\lim_{s\to0}E \biggl( \sup_{t\in[0,s]} \biggl| \frac{N^\Lambda
(t)}{v_\psi(t)} - 1 \biggr|^d \biggr)=0.
\end{equation}
\end{theorem}

The following consequence of Theorem \ref{Tsmalltime}
says that, among all the $\Lambda$-coalescents
such that $\Lambda[0,1]=1$,
Kingman's coalescent is extremal for
the speed of coming down from infinity.
\begin{coro}
\label{C:Kfast}
Assume (\ref{unit_mass}). Then with probability 1,
for any $\varepsilon>0$, and for all $t$ sufficiently small,
\[
N^\Lambda(t) \ge\frac2t (1-\varepsilon).
\]
\end{coro}
\begin{pf}
Without loss of generality assume that the $\Lambda$-coalescent comes down
from infinity.
To see how the corollary follows from Theorem \ref{Tsmalltime},
observe that since $e^{-qx}\le1-qx+ q^2x^2/2$ for
$x>0$,
%
%
\begin{equation}
\label{psi -g} \psi(q) \le\frac{q^2}{2}\int_{[0,1]} x^2\nu(dx)
\le
\frac{q^2}{2} \qquad\mbox{[due to (\ref{unit_mass})]}.
\end{equation}
Hence
%
%
\begin{equation}\label{E v nonint}
u_{\psi}(s)\ge\int_{s}^{\infty} \frac{2}{q^2} \,dq=\frac{2}s
\quad\mbox{and}\quad v_\psi(t)\geq\frac{2}t.
\end{equation}
Due to Theorem \ref{Tsmalltime},
$N^\Lambda(t)\sim v_\psi(t)$ as $t\to0$, implying that $N^\Lambda
(t)\ge
2(1-\varepsilon)/t$ with high probability for
all $t$ small.
\end{pf}
\begin{remark}
It is interesting to compare the last result with
the following fact shown in Angel et al. \cite{abhl}:
%
%
\begin{equation}\label{E:abhl}
\int_0^1 N^\Lambda(t) \,dt = \infty,
\end{equation}
regardless of the choice of the finite measure $\Lambda$.
Corollary \ref{C:Kfast} may be used to give an alternative proof of
(\ref{E:abhl}).
\end{remark}

The following result is interesting from the perspective of
applications in population genetics. More specifically, the total
length of the coalescent tree is relevant for predicting the number
of mutations in a large but finite sample. Assume (\ref{EequivaleG}), so that the coalescent comes down from infinity.
 Let $N^{\Lambda,n}$ denote
the number of blocks process of the restriction $\Pi^{(n)}$ with
initial state $\Pi_0^{(n)}=\{\{1\},\ldots,\{n\}\}$ as defined
at the beginning of Section \ref{S:Lacoal}. Let $\tau_n:=\inf\{s
>0\dvtx
N^{\Lambda}(s)\leq
n\}$, and let $H_n:=\{N^{\Lambda}(\tau_n)=n\}$ be the event that the
(unrestricted) $\Lambda$-coalescent ever attains a configuration with
exactly $n$ blocks. Then, due to the strong Markov property, the
conditional law of $(N^{\Lambda}(s+\tau_n), s\geq0)$ given
$\mathcal{F}_{\tau_n}$ on the event $H_n$, equals the law of
$N^{\Lambda,n}$.
Let $t_n=u_\psi(n)$ so that $v_\psi(t_n)=n$.
%
\begin{theorem}
For each $s>0$ we have \label{C:extraasym}
\[
\lim_{n\to\infty} \frac{\int_0^s N^{\Lambda,n}(t) \,dt}{\int_0^s
v_\psi(t_n+t) \,dt} =
\lim_{n\to\infty} \frac{\int_0^s N^{\Lambda,n}(t) \,dt}{\int_0^s E
(N^{\Lambda,n}(t)) \,dt} =1\qquad\mbox{in probability}.
\]
For Kingman and Beta coalescents [i.e., when $\Lambda$ is of the form
$\Lambda=\delta_0$ or (\ref{D:beta}) with $1< \alpha<2$], the above
convergence holds almost surely.
\end{theorem}

Let $\tau_1^n = \inf\{ t \ge0 \dvtx N^{\Lambda, n} (t) =1\}$,
so that $\int_0^{\tau_1^n} N^{\Lambda,n}(t) \,dt$ equals
the total length of the ($\Lambda$-)coalescent tree with
$n$ leaves.
Moreover, for any fixed $s>0$,
\[
\int_s^{\tau_1^n} N^{\Lambda,n}(t) \,dt \to\int_s^{\tau_1}
N^{\Lambda
}(t)\,
dt \qquad\mbox{almost surely}
\]
(see Section \ref{S:proof of theorem 5})
where the limit is a finite random variable.
Hence the above theorem yields the asymptotics for the total length of
the coalescent (genealogical) tree.
Some more detailed analysis is postponed until \cite{bbl2}.

Whereas Theorem \ref{Tsmalltime} is a law of large numbers-type
result for $N^\Lambda$, Theorem \ref{C:extraasym} is a law of large
numbers-type result for $\int_0^{\tau_1^n} N^{\Lambda,n}(t) \,dt$.
A central limit theorem for lengths of partial coalescent trees
is obtained by Delmas, Dhersin and Siri-Jegousse \cite{ddsj} (see
also \cite{sj}) for the Beta-coalescent case, similar questions for
general $\Lambda$-coalescents remain open.

\section{Martingale based arguments}
\label{S:martingale}
We now proceed toward the proof of Theorem \ref{Tsmalltime}.
The following easy-to-check facts will be used
in our analysis.
\begin{lemma}
\label{Lfewfacts}
The function $\psi\dvtx[0,\infty) \to\mathbb{R}_+$ of (\ref{D:psi}) is
(strictly) increasing on $[0,\infty)$, and convex on $(0,\infty)$.
Furthermore, for $v_\psi$, as in (\ref{Ev}), we have
$v_\psi'(s)=-\psi(v_\psi(s))$,
so that $v_\psi$ is decreasing with its derivative decreasing in
absolute value.
\end{lemma}

Due to Lemma \ref{Lsupport}, postponed until the next section,
we can, and will, suppose without loss of generality that
$\operatorname{supp}(\Lambda)\subset[0,1/4]$.
This assumption simplifies some technical estimates.

In this section we write $N$ instead of $N^\Lambda$ whenever
not in risk of confusion, and we also abbreviate $v=v_\psi$.
We start by observing that the function $v$ is the unique solution of
the following
integral equation:
%
%
\begin{equation}\label{IE for v}
\log(v(t)) - \log(v(z)) + \int_z^t \frac{\psi(v(r))}{v(r)} \,dr =
0\qquad
\forall0<z<t,
\end{equation}
with the ``initial condition'' $v(0+)=\infty$ [see Lemma \ref
{Lincreasing} for
properties of $\psi(q)/q$].
It is then natural to
consider, for each fixed $z>0$, the process
%
%
\begin{equation}
\label{EMart_almost} M(t):=\log(N(t)) - \log(N(z)) + \int_z^t
\frac{\psi(N(r))}{N(r)} \,dr,\qquad t \geq z.
\end{equation}
Let $n_0\ge1$ be fixed. Define
%
%
\begin{equation}
\label{E_taun0}
\tau_{n_0} :=\inf\{s>0\dvtx N(s) \leq n_0\}.
\end{equation}

The following proposition tells us that $M(t\wedge\tau_{n_0})$ is
``almost'' (up to a bounded drift correction, and integrability condition)
a martingale, with respect to
the natural filtration $(\mathcal{F}_t,t\ge0)$
generated by the underlying $\Lambda$-coalescent
process.
Its proof uses some general facts about
binomial distributions, with precise statements and arguments
postponed until the \hyperref[ap:binom]{Appendix}.
In particular, in the rest of this section
the parameter $n_0$
is taken to be the integer $n_0$ from Lemma \ref{Lbinlogcalc}.

As usual, $E[d X_s|\mathcal{F}_s]$ denotes the infinitesimal drift of
a continuous-time
process $(X_s,s\ge0)$ with respect to the filtration $\mathcal{F}$ at
time $s$.
Similarly, we denote by $E[(d X_s)^2|\mathcal{F}_s]$ the corresponding
infinitesimal
second moment. That is,
\[
\frac{E[d X_s|\mathcal{F}_s]}{ds}:= \lim_{\varepsilon\to0} \frac
1\varepsilon
E[X_{s+\varepsilon} - X_s | \mathcal{F}_s]
\]
and
\[
\frac{E[(d X_s)^2|\mathcal{F}_s]}{ds}:= \lim_{\varepsilon\to0}
\frac1\varepsilon
E[(X_{s+\varepsilon} - X_s)^2 | \mathcal{F}_s].
\]
\begin{proposition} \label{martingale_estimates}
There exists some deterministic $C<\infty$ such that
%
%
\begin{equation}\label{Eha}
E[d\log(N(s))|\mathcal{F}_s]= \biggl(-\frac{\psi(N(s))}{N(s)} +
h(s)
\biggr) \,ds,
\end{equation}
where $(h(s), s\ge z)$ is an $\mathcal{F}$-adapted process such that
${\sup_{s\in[z,z \wedge\tau_{n_0}]}} |h(s)|\le C$, and
\[
E[[d\log(N(s))]^2|\mathcal{F}_s] \mathbf{1}_{\{s \leq\tau_{n_0}\}}
\leq C \,ds\qquad
\mbox{almost surely.}
\]
Both estimates are valid uniformly over $z>0$.
\end{proposition}
\begin{pf}
To prove the proposition, it suffices to
show that for each $s>0$, we have on $\{N(s) \geq n_0\}$,
%
%
\begin{equation}\label{eqdriftvar}
\biggl|\frac{E(d \log(N(s))|\mathcal{F}_s)}{ds} + \frac{\psi(N(s))}{N(s)}
\biggr|=
|h(s)|= O \biggl(\int_{[0,1/4]}p^2 \nu(dp) \biggr)
\end{equation}
and
%
%
\begin{equation}\label{eqdriftvar2}
E([d \log(N(s))]^2|\mathcal{F}_s) =O \biggl(\int_{[0,1/4]}p^2 \nu
(dp) \biggr)\,ds,
\end{equation}
where $O(\cdot)$ can be taken uniformly in $s$.
Note that the finite integrals above are in fact taken over $[0,1]$,
since $\nu(dx)=\nu(dx) \mathbf{1}_{\{x\in[0,1/4]\}}$ by assumption.

Recall the Poisson point process construction of
Section \ref{S:Lacoal} and fix $n\geq n_0$.
If $\Lambda(\{0\})=0$, then on the event $\{N(s)
=n\}$ an atom of size $p$ arrives at rate $\nu(dp) \,ds$, and given
that, $\log{N(s)}=\log{n}$ jumps to $\log(B_{n,p}+\mathbf{1}_{\{
B_{n,p} <n\}})$ where $B_{n,p}$ has $\operatorname{Binomial}(n,1-p)$ distribution. Hence
we have
\[
E(d \log(N(s))|\mathcal{F}_s)=\int_{[0,1]} E \biggl[\log\frac
{B_{n,p}+\mathbf{1}_{\{B_{n,p} <n\}}}{n} \biggr] \nu(dp) \,ds.
\]
In the general case where $\Lambda(\{0\})=c\in(0,1)$, we have
on the same event
\begin{eqnarray*}
\frac{E(d \log(N(s))|\mathcal{F}_s)}{ds} &=& (1-c) \int_{[0,1]}
E \biggl[\log\frac{B_{n,p}+\mathbf{1}_{\{B_{n,p} <n\}}}{n} \biggr]
\nu(dp)\\
&&{}
+c\pmatrix{n \cr2} \log\frac{n-1}{n}.
\end{eqnarray*}
Let $\psi_0(q)=q^2/2$ be the function $\psi$ corresponding to the
atomic $\Lambda(dx)=\delta_0(dx)$.
Note that
\[
c\pmatrix{n \cr2} \log\frac{n-1}{n} = - \frac{c}{2} n + \frac{c}{4} +
O(1/n) =
-c\frac{\psi_0(n)}{n} + \frac{c}{4} + O(1/n).
\]
In view of (\ref{Epsi_atom}), the estimate (\ref{eqdriftvar}) will
follow by Lemma \ref{Lbinlogcalc} in the \hyperref[ap:binom]{Appendix}
provided that
%
%
\begin{equation}
\label{Eexppowbd}\quad
|np-1+(1-p)^n - (e^{-np}-1 + np)| = |(1-p)^n -e^{-np}|\leq C n p^2
\end{equation}
for all $n\geq n_0$, $p\leq1/4$ and for some $C<\infty$. Note that
$e^{-np}>(1-p)^n$, and, in fact,
\[
e^{-np} - (1-p)^n = e^{-np} \biggl(1- \exp\biggl\{-n \biggl(\frac
{p^2}{2}+\frac{p^3}{3}+\cdots\biggr) \biggr\} \biggr).
\]
Therefore, for $p\leq1/4$
we have
\begin{eqnarray*}
1- \exp\biggl\{-n \biggl(\frac{p^2}{2}+\frac{p^3}{3}+\cdots
\biggr) \biggr\} &\leq& 1-\exp\biggl\{-\frac{n}{2} (p^2 + p^3 +
\cdots) \biggr\} \\
& \le& 1 - \exp\biggl(-\frac{2}{3} np^2 \biggr) \\
& \le& \frac{2}{3} np^2;
\end{eqnarray*}
%
hence both (\ref{Eexppowbd}) and (\ref{eqdriftvar}) hold.

To bound the infinitesimal variance on the event $\{N(s) =n\}$,
use the second estimate in Lemma \ref{Lbinlogcalc}, together with the fact
\begin{eqnarray*}
\frac{E([d \log(N(s))]^2|\mathcal{F}_s)}{ds} &\leq&
(1-c) \int_{[0,1]} E \biggl[\log^2\frac{B_{n,p}+\mathbf{1}_{\{
B_{n,p} <n\} }}{n} \biggr] \nu(dp)\\
&&{}
+ O \biggl(\frac{1}{n^2} \biggr) \pmatrix{n \cr2}.
\end{eqnarray*}
Finally, note that both bounds (\ref{eqdriftvar}) and (\ref
{eqdriftvar2}) are uniform in the choice of $z$.
\end{pf}

\subsection[Proof of Theorem 1]{Proof of Theorem \protect\ref{Tsmalltime}}
\label{S:smalltime}
Recall the process $M$ from (\ref{EMart_almost}) and define
%
%
\begin{equation}
\label{Elocmart}
M_z'(t)\equiv M'(t):= M(t\wedge\tau_{n_0}) + \int_z^{t\wedge\tau
_{n_0}} h(r)\,dr,\qquad
t\geq z,
\end{equation}
so that $M_z'$ has martingale increments due to Proposition \ref
{martingale_estimates}.
A general property of the Doob--Meyer martingale correction (that one
can check easily)
implies that
$E([dM_z'(t)]^2|\mathcal{F}_t)\leq E[[d\log(N(t))]^2|\mathcal
{F}_t]$, so that
%
%
\begin{equation}\label{E varM estim}
E\bigl(M_z'(s) - M_z'(z)\bigr)^2 \leq C (s-z)\qquad\forall0<z<s,
\end{equation}
where $C$ is the constant from Proposition \ref{martingale_estimates}.

Define a family of deterministic functions $(v_x, x\in\mathbb{R})$ by
\[
v_x(t)= v(t+x),\qquad t\geq-x,
\]
and note that each $v_x$ satisfies an appropriate analogue of (\ref{IE
for v}) on its entire domain,
namely, $ v_x(-x+)=\infty$ and
%
%
\begin{equation}\label{IE for v shift}
\log(v_x(t)) - \log(v_x(z)) + \int_z^t \frac{\psi
(v_x(r))}{v_x(r)} \,dr = 0\qquad
\forall{-x}<z<t.
\end{equation}
For each fixed $z>0$ and each $x>-z$, define
\begin{eqnarray*}
M_{z,x}(t)&:=&
\log{\frac{N(t)}{v_x(t)}} -
\log{\frac{N(z)}{v_x(z)}}\\
&&{}+
\int_z^t
\biggl[
\frac{\psi(N(r))}{N(r)} -
\frac{\psi(v_x(r))}{v_x(r)}
+h(r) \biggr] \,dr,\qquad t \ge z,
\end{eqnarray*}
where $h$ is given in (\ref{Eha}).

Moreover, given $X\in\mathcal{F}_z$ such that $P(X>-z)=1$, we can consider
the process $M_{z,X}$.
The advantage of this approach will be apparent soon.

For fixed $z>0$, the processes $M_z',M_{z,x}$ and $M_{z,X}$ are all
adapted to the filtration
$(\mathcal{F}_r, r\geq z)$.
\begin{remark}
Strictly speaking, the processes $M_z',M_{z,x}$ and $M_{z,X}$ defined
above are \textit{local} martingales
(see \cite{revuzyor} Chapter II or \cite{rogerswilliams}, Chapters VI, 31--34
for definition and first properties)
since we do not know a priori whether $\log(N(t))$ has finite expectation.
However, the optional stopping and Doob moment estimates that we apply below
still hold in this more general setting.
\end{remark}
\begin{lemma}
\label{Lincreasing}
The function $q\mapsto\psi(q)/q$ is increasing.
\end{lemma}
\begin{pf}
Note that $q \mapsto\psi(q)/q$ is smooth, and that its derivative at
$q$ equals
\begin{eqnarray*}
\frac{\psi'(q) q - \psi(q)}{q^2} &=&
\frac{\int(1 - (xq +1) e^{-qx}) \nu(dx)}{q^2}\\
&=&
\frac{\int(1 - (xq +1) e^{-qx})/x^2 \Lambda(dx)}{q^2}.
\end{eqnarray*}
It is a simple matter to check that the integrand in the numerator
is positive for all $x>0$, and that its limit as $x\to0$ is $q^2/2$,
so again it is positive.

The reader is invited to verify in a similar manner
that $\lim_{q\to\infty} (\psi(q)/q)''=-\int_{[0,1]} e^{-qx}
x\Lambda
(dx)$ which implies that
$q\mapsto\psi(q)/q$ is asymptotically concave.
Our argument does not make use of this fact.
\end{pf}

The following deterministic lemma is a crucial step in our analysis. It
overcomes the need
for a priori estimates necessary for the method of \cite{darnor} to
apply, as discussed
in the \hyperref[sec1]{Introduction}.
\begin{lemma}
\label{Ldeterministic}
Suppose $f,g\dvtx[a,b]\mapsto{\mathbb R}$ are deterministic c\`{a}dl\`{a}g
functions
such that
%
%
\begin{equation}
\label{Ecoupledfun}
{\sup_{x\in[a,b]}} \biggl|f(x) + \int_a^x
g(u) \,du \biggr| \leq c
\end{equation}
for some $c<\infty$.
If, in addition,
$f(x)g(x) > 0$, $x\in[a,b]$ whenever $f(x) \neq0$,
then both
\[
\sup_{x\in[a,b]} \biggl|\int_a^x g(u) \,du \biggr|\leq c \quad\mbox
{and}\quad
{\sup_{x\in[a,b]} }|f(x)| \leq2c.
\]
%
\end{lemma}
\begin{pf}
Due to the assumptions, we know that at any point $x$ where $f(x)$
is positive (resp. negative) $h(x):=\int_a^x g(u) \,du$ is
increasing (resp. decreasing) from the right. Define $t_1
:=\min\{x \in[a,b]\dvtx|h(x)| > c\}$, with the convention that
$t_1=b$ if this set is empty. Suppose $t_1<b$. By
continuity of $h$, it must be that $|h(t_1)|=c$.
Without loss of generality, assume
%
%
\begin{equation}
\label{Eassumepos}
h(t_1)=c\quad\mbox{and hence}\quad h(t_1+\varepsilon)>c
\end{equation}
for all small enough $\varepsilon>0$. Having $f(t)<0$ for all $t\in
(t_1,t_1+\varepsilon)$ would imply that $h$ is decreasing on that same
interval, contradicting (\ref{Eassumepos}). Therefore, there exists
$t\in(t_1,t_1+\varepsilon)$ such that $f(t)\ge0$. But since $h(t)>c$ by
(\ref{Eassumepos}), this would in turn contradict
(\ref{Ecoupledfun}). Hence it must be $t_1=b$, so that the uniform
bound on $|h|$ holds, which together with (\ref{Ecoupledfun})
implies the uniform bound on $|f|$.
\end{pf}

Since
$N(t) \to\infty$ as $t\to0$, almost surely, we have
\[
P(\tau_{n_0}>0)=1\quad\mbox{or equivalently}\quad\lim_{s\to0}P(\tau
_{n_0} \leq s) = 0.
\]
Therefore, for any family $(Y_s, s>0)$ of random variables, we
have\break
$\lim_{s\to0, s\leq\tau_{n_0}} Y_s= \lim_{s\to0} Y_s$, almost
surely (in the sense that whenever one of the limits exists so does
the other). Without loss of generality we will henceforth write
$M_{z,x}(t)$ instead of $M_{z,x}(t\wedge\tau_{n_0})$, $t\in[z,s]$
instead of $t\in[z,s\wedge\tau_{n_0}]$, etc.

Fix $\alpha^*\in(0,1/2)$.
By Doob's
$L^2$-inequality for martingales and (\ref{E varM estim}) 
we have
%
%
\begin{eqnarray}
\label{Eboundprime}\quad
P\Bigl({\sup_{t\in[z,s]} }|M_z'(t)-M_z'(z)| > s^{\alpha^*}\Bigr) &\le& s^{-2
\alpha^*}
\sup_{t \in[z,s]}
E\bigl[\bigl(M_z'(t)-M_z'(z)\bigr)^2\bigr]\nonumber\\[-8pt]\\[-8pt]
&\leq& s^{-2
\alpha^*} C(s-z)= O(s^{1-2\alpha^*}).\nonumber
\end{eqnarray}
Denote by
\[
A_z'(s)\equiv A_z':= \Bigl\{{\sup_{t\in[z,s]}} |M_z'(t)-M_z'(z)| \leq
s^{\alpha^*}\Bigr\}
\]
the complement of the above event. Henceforth we assume that $s <
(1/C)^{1/({1-\alpha^*})}$. Note that then $\int_z^s h(r) \,dr \leq
\int_z^s C \,dr \le Cs \leq s^{\alpha^*}$. So we obtain that on
$A_z'$ [hence with probability greater than $1- O(s^{1-2\alpha^*})$],
\[
\sup_{t\in[z,s]} \biggl|\log{N(t)} - \log{N(z)} + \int_z^t
\frac{\psi(N(r))}{N(r)} \,dr \biggr| \leq2 s^{\alpha^*}.
\]
We conclude that $A_z'\subset A_z$, where
%
%
\begin{eqnarray}
\label{def_Az}\quad
A_z(s)&\equiv& A_z\nonumber\\[-8pt]\\[-8pt]
:\!&=&
\biggl\{\sup_{t_1,t_2\in[z,s]}
\biggl|\log{N(t_2)} - \log{N(t_1)} + \int_{t_1}^{t_2}
\frac{\psi(N(r))}{N(r)} \,dr \biggr| \leq4 s^{\alpha^*}
\biggr\}.\nonumber
\end{eqnarray}
The advantage of the new definition is that $A_{z_1}\subset A_{z_2}$
whenever $z_1\leq z_2\leq s$.
Moreover, the bound in (\ref{Eboundprime}) is uniform in $z\in(0,s)$, hence
the decreasing property of probability measures implies
\begin{eqnarray*}
P\biggl(\bigcap_{z\in(0,s)} A_z\biggr)
&=&
P \biggl(\sup_{t_1,t_2\in(0,s]}
\biggl|\log{N(t_2)} - \log{N(t_1)} + \int_{t_1}^{t_2}
\frac{\psi(N(r))}{N(r)} \,dr \biggr| \leq4 s^{\alpha^*} \biggr)
\\
&=& 1- O(s^{1-2\alpha^*}).
\end{eqnarray*}
Let $X_z$ be the random variable defined by
%
%
\begin{equation}
\label{def_Xz}
N(z)=v(X_z+z)= v_{X_z}(z).
\end{equation}
\begin{lemma}
\label{L conv Xz}
We have $\lim_{z\to0} X_z = 0$, almost surely.
\end{lemma}
\begin{pf}
Since $N$ is nonincreasing and $v$ is (strictly) decreasing, it is
easy to see that $(X_z+z, z > 0)$ is also a nondecreasing
process, almost surely. Therefore $\lim_{z\to0} X_z +z \geq0$
exists, almost surely. Moreover, the above limit equals $0$ with
probability $1$, since $X_z+z = u(N(z))$, and since
$P(N(0+)=\infty)=1$ and $\lim_{x\to\infty}u(x)=0$.
\end{pf}

Due to (\ref{IE for v shift}) and (\ref{def_Az}),
we have, in particular, that
\begin{eqnarray*}
&&A_z= \biggl\{\sup_{t_1,t_2\in[z,s]}
\biggl|\log{\frac{N(t_2)}{v_{X_z}(t_2)}} -
\log\frac{N(t_1)}{v_{X_z}(t_1)}\\
&&\hspace*{68.8pt}{} + \int_{t_1}^{t_2}
\biggl[\frac{\psi(N(r))}{N(r)} - \frac{\psi(v_{X_z}(r))}{v_{X_z}(r)}
\biggr] \,dr \biggr| \leq4 s^{\alpha^*} \biggr\}.
\end{eqnarray*}
After plugging in $t_1=z$, we obtain
\[
A_z \subset
\biggl\{\sup_{t\in[z,s]}
\biggl|\log{\frac{N(t)}{v_{X_z}(t)}} + \int_{z}^{t}
\biggl[\frac{\psi(N(r))}{N(r)} - \frac{\psi
(v_{X_z}(r))}{v_{X_z}(r)} \biggr] \,dr \biggr| \leq4 s^{\alpha^*}
\biggr\}.
\]
Lemma \ref{Lincreasing} implies the hypotheses of Lemma
\ref{Ldeterministic} omega-by-omega (with $a=z$, \mbox{$b=s$} and
the obvious choice of $f$ and $g$), therefore
%
%
\begin{equation}
\label{EAzbound}
A_z(s) =A_z \subset
\biggl\{
\sup_{t\in[z,s]} \biggl|\log{\frac{N(t)}{v_{X_z}(t)}} \biggr|
\leq8 s^{\alpha^*}
\biggr\}
.
\end{equation}
By fixing $t<s$ and varying $z\in(0,t]$
[note that
$\log\frac{v_{X_z}(t)}{v_{X_{z'}}(t)}=
\log\frac{N(t)}{v_{X_{z'}}(t)}-\log\frac{N(t)}{v_{X_z}(t)}$]
we obtain
\[
\bigcap_{z\in(0,s)} A_z(s) \subset\biggl\{
\sup_{z,z' \in(0,s), t\in[z\vee z',s]}
\biggl|\log{\frac{v_{X_z}(t)}{v_{X_{z'}}(t)}} \biggr|
\leq16 s^{\alpha^*}
\biggr\},
\]
which together with (\ref{EAzbound}) implies
\[
\bigcap_{z\in(0,s)} A_z(s) \subset\biggl\{
\sup_{t\in(0,s]}
\biggl|\log{\frac{N(t)}{\lim_{n}v_{X_{z_n}}(t)}} \biggr|
\leq24 s^{\alpha^*}
\biggr\},
\]
where $(z_n)_{n\geq1}$ is any given
deterministic sequence of strictly
positive\break numbers converging to $0$.
Due to Lemma \ref{L conv Xz},
the continuity of $v$ implies
$\lim_{n\to\infty}v_{X_{z_n}}(t)=v(t)$, $\forall t\in(0,s]$, almost surely.
To summarize, we have just proved:
\begin{proposition}
\label{Psummar}
If $\operatorname{supp}(\Lambda)\subset[0,1/4]$, then
\[
P \biggl(
\sup_{t\in(0,s \wedge\tau_{n_0}]}
\biggl|\log\frac{N(t)}{v(t)} \biggr|
\leq24 s^{\alpha^*}
\biggr)\geq P\biggl(\bigcap_{z\in(0,s)} A_z(s)\biggr)= 1- O(s^{1-2\alpha^*}).
\]
\end{proposition}

Theorem \ref{Tsmalltime}
now follows
due to the Borel--Cantelli lemma, after choosing a deterministic
sequence $(s_m)_{m\geq1}$
of strictly positive numbers converging to $0$ sufficiently fast so
that $\sum_m (s_m)^{1-2\alpha^*}<\infty$.
\begin{remark}
The fixed scale assumption $\Lambda[0,1] (= \Lambda[0,1/4]) = 1$ has
not been
used in the above argument.
\end{remark}

\subsection{Relaxing assumptions on $\operatorname{supp}(\Lambda)$}
\label{S:relax}
Given a probability measure $\Lambda$ on $[0,1]$ and a positive $\eta
\leq
1$, define
its restriction $\Lambda_\eta$ by
\[
\Lambda_\eta(dx) = \Lambda(dx)\mathbf{1}_{[0,\eta]}(dx).
\]
For each $\eta\in(0,1]$, denote by $\psi_\eta$ the function $\psi
_{\Lambda_\eta}$
that corresponds to $\Lambda_\eta$ [cf. (\ref{D:psi})], and by
$v_\eta$
the corresponding rate
function from (\ref{Ev}).
\begin{lemma}
\label{Lsupport}
All the $\Lambda_\eta$-coalescents, where $\eta\in(0,1]$, have the same
speed of CDI. Moreover, for any fixed $\eta\in(0,1)$,
%
%
\begin{equation}
\label{Evveta}
\lim_{t \to0} \frac{v(t)}{v_\eta(t)}=1,
\end{equation}
so it suffices to prove Theorem \ref{Tsmalltime} for
one $\eta\in(0,1)$ in order to prove it for all $\eta\in(0,1]$.
\end{lemma}
\begin{pf}
Fix $\eta\in(0,1)$.
Assume first that $\Lambda(\{0\})=0$.
Then it is easy to see that one can find a coupling of the two
coalescent processes defined by
$\Lambda$ and by $\Lambda_\eta$, respectively,
such that the corresponding coalescent block counting
processes $N^\Lambda$ and $N^{\Lambda_\eta}$
coincide for all $t\in(0,T_\eta)$
where \mbox{$P(T_\eta>0) =1$}.
Namely, recall the PPP construction of Section \ref{S:Lacoal}
and set
$
T_\eta:=\min\{t>0\dvtx(t,p) \mbox{ is an atom of }\pi\mbox{ and
}p>\eta\}.
$

If $\Lambda(\{0\})>0$, let $\Lambda'(dx)=\Lambda(dx)\mathbf{1}_{(0,1)}(x)$, and
note that the PPP-based construction of $\Lambda'$-coalescent
can be enriched by superimposing pairwise coalescent
events at rate $\Lambda(\{0\})$
thus yielding a construction of $\Lambda$-coalescent.
Again, one can couple such constructions of
$\Lambda$-coalescent and $\Lambda_\eta$-coalescent so that the two
processes agree
until $T_\eta$ as discussed above.

To prove the lemma, it now suffices to show (\ref{Evveta}) for any
fixed $\eta\in(0,1)$.
Note that
we trivially have $v(t)\leq v_\eta(t)$ for all $t>0$,
since $\psi_\eta(q)\leq\psi(q)$ for all $q>0$.
Moreover, 
\[
\psi_\eta(q)=\psi(q)- a_\eta q + b_\eta+ O(e^{-q \eta}),
\]
where $a_\eta:= \int_{(\eta,1]} ({1}/{x}) \Lambda(dx)$ and $b_\eta
:=\int_{(\eta,1]} ({1}/{x^2}) \Lambda(dx)$.
Therefore, for any $0\leq z\leq t$,
\[
\log{\frac{v(t)}{v_\eta(t)}} -
\log{\frac{v(z)}{v_\eta(z)}}+
\int_z^t
\biggl[
\frac{\psi(v(r))}{v(r)} -
\frac{\psi(v_\eta(r))}{v_\eta(r)}
+h_z(r) \biggr] \,dr = 0,
\]
where $h_z(r)$ is now a deterministic function, bounded by a fixed
constant $C$, uniformly over $z$.
The rest of the argument is a deterministic (and easier) analogue of
the argument
given in Section \ref{S:smalltime}.
We leave it to an interested reader.
\end{pf}

If $\Lambda(\{0\})>0$, then the size of the atom at $0$
determines the speed of CDI. More precisely, we have:
\begin{coro}
If $\Lambda(\{0\})=c>0$, then for all $\eta\in(0,1]$,
\[
v_\eta(t) \sim\frac{2}{c t} ,\qquad t \to0.
\]
\end{coro}
\begin{pf}
Denote by $v_0$ the above function $2/(c t)$ and note that
it corresponds to $\Lambda(dx)= c\delta_0(dx)$ and $\psi_0(q)= \frac
{c q^2}{2}$,
in terms of (\ref{Ev}).
Next note that if $\eta\in(0,1]$, then
\[
\psi_\eta(q)=\frac{c q^2}{2}+ f(q)= \psi_0(q) + f(q),
\]
where $f(q)= o(q^2)$ is a nonnegative function.
In particular, $v_\eta(t)\leq v_0(t)$, $t >0$.
Moreover, since for any $\varepsilon$, we can find $q(\varepsilon
)<\infty$, such that
\[
\psi_\eta(q)\leq\frac{c (1+\varepsilon) q^2}{2}\qquad \mbox{for all }
q \geq
q(\varepsilon).
\]
We have by the same reasoning,
$v_\eta(t)\geq v_0(t)/(1+\varepsilon)$ for all sufficiently small $t$.
Letting $\varepsilon\to0$ implies the statement.
\end{pf}

\subsection[Proof of Theorem 2]{Proof of Theorem \protect\ref{Tsmalltime-mom}}
Assume that the parameter $n_0$ is the maximum of the
corresponding quantities from Lemmas \ref{Lbinlogcalc} and \ref
{Lbinlogcalc1}.
Assume initially that $\operatorname{supp}(\Lambda) \subset[0,1/4]$ and
fix $z>0$.
With the notation of Section \ref{S:smalltime} in mind,
let $M_{z,X_z}\equiv M$
be the process given by
%
%
\begin{eqnarray}\label{E M thm2}
M_t &:=& \log\frac{N(t \wedge\tau_{n_0})}{v(X_z+t \wedge\tau_{n_0})}
\nonumber\\[-8pt]\\[-8pt]
&&{}
+ \int_z^{t\wedge\tau_{n_0}} \biggl(\frac{\psi(N(r))}{N(r)} -
\frac{\psi(v(X_z+r))}{v(X_z+r)} +h(r) \biggr) \,dr,\qquad t \geq
z.\nonumber
\end{eqnarray}
Then $M_z=0$, and due to Proposition \ref{martingale_estimates}, $M$
is a martingale (in the sense that $M_t$ is an integrable random
variable, $t\geq z$). Note that here we use $M$ as abbreviation; the
above process should not be confounded with $M$ from (\ref{EMart_almost}).

We next obtain better estimates on the tails of the distribution of
$M_t$, via
an analogue of Hoeffding's inequality \cite{hoeffding} for discrete
martingale sums.
Since $M$ has only downward jumps, a simple case of a general result
of Barlow, Jacka and Yor (\cite{barlowetal}, Proposition 4.2.1; see
also \cite
{delapena})
implies that for any $c>0$,
\[
S^{(c)}:= \biggl(\exp\biggl\{ c M_t - \frac{c^2 C (t-z)}{2} \biggr\},
t\geq z \biggr),
\]
is a supermartingale started from $S_z^{(c)}=1$,
with respect to the usual filtration $\mathcal{F}$.
Note that $D_t$ in
\cite{barlowetal,delapena} corresponds to $E[(dM_t)^2| \mathcal{F}_t
]$ in
our notation,
and that
$C$ is the uniform upper bound from Proposition \ref{martingale_estimates}.

Fix some $x \in\mathbb{R}_+$.
Let $c=x/(C (s-z))$, and $y=\exp\{cx/2\}= \exp\{ cx -c^2 C (s-z)/2)\}$,
and let $T_y = \inf\{ t \ge z\dvtx S^{(c)}_t > y\}$. Since $S^{(c)}$ only
has downward jumps,
it must be $S^{(c)}_{T_y} = y$ on $\{T_y<\infty\}$.
Since $S^{(c)}$ is supermartingale, using optional stopping at
$T_y \wedge s$, we have
\begin{eqnarray*}
1=E \bigl(S^{(c)}_z \bigr) &\ge& E \bigl(S^{(c)}_{T_y\wedge s} \bigr)
\\
& = & y P(T_y \le s) + E \bigl(S^{(c)}_s \mathbf{1}_{T_y > s} \bigr)\\
&\ge & yP(T_y \le s).
\end{eqnarray*}
It follows that
\begin{eqnarray*}
P \Bigl(\sup_{t\in[z,s]} M_t > x \Bigr) &\le& P \Bigl(\sup_{t\in
[z,s]} S_t^{(c)} > e^{cx- c^2 C (s-z)/2} \Bigr)\\
& \le & P(T_y \le s) \\
& \le &\frac1y = \exp\biggl\{-\frac{x^2}{2C(s-z)} \biggr\}.
\end{eqnarray*}

In order to obtain the ``left tails'' we use
\cite{barlowetal} Proposition 4.2.1 in a less trivial sense.
If $c>0$, then
\[
S^{(-c)}:= \biggl(\exp\biggl\{ -c M_t - \frac{c^2 C (t-z)}{2} - \frac
{c^2}{2}\sum_{s\leq t} (\Delta_s M)^2 \biggr\}, t\geq z \biggr),
\]
is a supermartingale where $\Delta_s M = M(s)-M(s-)= \Delta_s \log
{N(s\wedge\tau_{n_0})}$.
Define
\begin{eqnarray*}
E^{(c)}(t)&:=& \exp\biggl\{c \sum_{t\in[z,s]} (\Delta_t M)^2
-e^{9c/4} K_0 (t-z) \biggr\}\\
&=& \exp\biggl\{c \sum_{t\in[z,s]} \bigl(\Delta_t \log{N(s\wedge
\tau_{n_0})}\bigr)^2 -e^{9c/4} K_0 (t-z) \biggr\},
\end{eqnarray*}
where $K_0$ is the constant from Lemma \ref{Lbinlogcalc1}.
Due to Lemma \ref{Lbinlogcalc1},
we have that for each $c>0$, the process
$
(E^{(c)}(t), t\geq z)
$
is a nonnegative super-martingale started from $E^{(c)}(z)=1$.
Indeed, it is easy to verify in the sense of calculations
of Proposition~\ref{martingale_estimates} that
\begin{eqnarray*}
&&E\bigl(d E^{(c)}(t)| \mathcal{F}_t\bigr) \\
&&\qquad= E^{(c)}(t)\cdot
E[\exp\{c(\Delta_t M)^2\}-1| \mathcal{F}_t ]-e^{9c/4} K_0 \cdot E^{(c)}(t)
\,dt\\
&&\qquad\leq E^{(c)}(t)\cdot\biggl[\sum_{n\geq n_0}
\mathbf{1}_{\{N(t)=n\}}\int_{[0,1/4]} (e^{9c/4} K_0 p^2) /p^2
\Lambda(dp)
-e^{9c/4} K_0 \biggr] \,dt\\
&&\qquad = 0,
\end{eqnarray*}
almost surely.
To include the case $\Lambda(\{0\})>0$
in the above calculation, note that by a
standard estimate (\ref{Ecalc_fact}) and Taylor's series expansion,
\[
\pmatrix{n \cr2} \bigl(\exp\bigl\{ c \log^2\bigl((n-1)/n\bigr) \bigr\} -1\bigr) =\frac{c}{2} +O
\biggl(\frac{c}{n} + \frac{e^c}{n^2} \biggr).
\]
Without loss of generality one can assume that both $K_0\geq1$ and
$c/2 + O ({c}/{n} + {e^c}/{n^2} )\leq e^{9c/4}$ for $n\geq
n_0$ and all $c >0$.

Then for $x>0$, we have
\begin{eqnarray*}
&&P \Bigl(\inf_{t\in[z,s]} M_t <- x
\Bigr)\\
&&\qquad\leq
P \biggl(\inf_{t\in[z,s]} M_t < -x,
c^2 \sum_{t\in[z,s]} (\Delta_s M)^2 \leq c x \biggr)\\
&&\qquad\quad{} +
P \biggl(\sum_{t\in[z,s]} (\Delta_t M)^2 > x/c \biggr)\\
&&\qquad\leq P \Bigl(\sup_{t\in[z,s]} S_t^{(-c)} > e^{cx/2- c^2 C
(s-z)/2} \Bigr)
\\
&&\qquad\quad{}
+P \Bigl(\sup_{t\in[z,s]} E^{(c^2)}(t)> e^{x c - e^{9c^2/4} K_0
(s-z)} \Bigr)\\
&&\qquad\leq e^{-cx/2 +c^2 C (s-z)/2} + e^{-x c + e^{9c^2/4} K_0
(s-z)}.
\end{eqnarray*}
We plug in $c=\frac23 \sqrt{\log[x/(K_0 (s-z))]}$ [here we assume
that $x
> 2 K_0 (s-z)$]. Since in each exponent the second term is negligible
when compared
to the first, we get the sub-exponential estimate
\[
P \Bigl(\inf_{t\in[z,s]} M_t <- x \Bigr) = O\bigl (r(x;s-z) \bigr),
\]
where
\[
r(x;s):= \exp\bigl\{-x \sqrt{\log[x/(K_0 s)]}/4 \bigr\}.
\]
Now another omega-by-omega application of
Lemmas \ref{Lincreasing} and \ref{Ldeterministic} yields
\begin{eqnarray*}
1-O\bigl(r(x;s-z)\bigr) &\leq& P \Bigl( \sup_{t\in[z,s]} |M_t| \leq x \Bigr)
\\
&\leq&P \biggl(\sup_{t\in[z,s]} \biggl|\log\frac{N(t\wedge\tau
_{n_0})}{v(X_z+t\wedge\tau_{n_0})}
\biggr|\leq2(x + C s)
\biggr).
\end{eqnarray*}
Since $\lim_{z \to0}v(X_z+t) = v(t)$
as argued before, in the limit we obtain
%
%
\begin{equation}
\label{E exp est}
P \biggl(\sup_{t\in[0,s]} \biggl|\log\frac{N(t\wedge\tau
_{n_0})}{v(t\wedge\tau_{n_0})} \biggr| \leq2(x + C s) \biggr) \geq
1- O(r(x;s)).
\end{equation}
%
Note that since $N$ is an integer-valued process and $v$ is a
decreasing function,
$\inf_{t\in[0,s]}\log(N(t)/v(t)) \geq\inf_{t\in[0,s\wedge\tau
_{n_0}]}\log(N(t)/v(t)) -\log{n_0}$,
almost surely.
Now (\ref{E exp est}) together with the observation $N(t)\leq
N(t\wedge\tau_{n_0})$
implies that the random variable
\[
\Xi_s:=\sup_{t\in[0,s]} \biggl| \log\frac{N(t)}{v(t)} \biggr|=
\log\biggl( \sup_{t\in[0,s]} \biggl|\frac{N(t)}{v(t)} \biggr|
\vee\sup_{t\in[0,s]} \biggl|\frac{v(t)}{N(t)} \biggr|
\biggr)
\]
satisfies
$P(\Xi_s >x) = O(r(x;s))$,
hence
\[
P \biggl( \sup_{t\in[0,s]} \biggl|\frac{N(t)}{v(t)} \biggr| \geq y
\biggr) \leq
O \biggl(\frac{1}{y^{\sqrt{\log\log(y)-\log(K_0s)}/4}} \biggr)
\qquad\mbox{as } y\to\infty.
\]
In particular,
for any $d\geq1$, we can find a constant $D(d)<\infty$
such that
%
%
\begin{equation}
\label{E all moments}
E \biggl(\sup_{t\in[0,s]} \biggl|\frac{N(t)}{v(t)} \biggr|^d \biggr)
< D(d),
\end{equation}
hence (for a possibly
different constant) $E ({\sup_{t\in[0,s]}} |N(t)/v(t)-1|^d) < D(d)$.
Now the almost sure convergence of Theorem \ref{Tsmalltime} combined with
an application of dominated convergence theorem completes the argument.

For the case of general $\operatorname{supp}(\Lambda)$, recall the
notation of Section
\ref{S:relax}.
In addition, denote by $N_{1/4}(t)$
the number of blocks process corresponding to $\Lambda_{1/4}$.
Due to the coupling construction used in the argument of Lemma \ref
{Lsupport}, we have
\[
N_{1/4}(t) \geq N(t),\qquad t\geq0,
\]
and moreover,
\[
\sup_{t\in[0,s]} \frac{v_{1/4}(t)}{v(t)} < \infty.
\]
Therefore estimate (\ref{E all moments}), established for the
$\Lambda_{1/4}$-coalescent,
will imply the same estimate
[with possibly different constant $D(d)$] for the $\Lambda$-coalescent.

\subsection[Proof of Theorem 5]{Proof of Theorem \protect\ref{C:extraasym}}\label
{S:proof of theorem
5}
Recall the notation $t_n=u_\psi(n)=u(n)$ introduced before the
statement of
Theorem \ref{C:extraasym}.
It suffices to show that any subsequence $(n_k)_{k\geq1}$ contains a
further subsequence
$(n_{k(j)})_{j\geq1}$ such that
%
%
\begin{eqnarray}
\label{Esuffic}\qquad
\lim_{j\to\infty} \frac{\int_0^s N^{\Lambda,n_{k(j)}}(t)
\,dt}{\int
_0^s v(t_{n_{k(j)}}+t) \, dt} = 1
 = \lim_{j\to\infty} \frac{\int_0^s N^{\Lambda,n_{k(j)}}(t)
\,dt}{\int
_0^s E (N^{\Lambda,n_{k(j)}}(t)) \,dt}\nonumber\\[-8pt]\\[-8pt]
\eqntext{\mbox{almost
surely}.}
\end{eqnarray}
For $t \geq0$, define
%
%
\begin{eqnarray}\label{E M thm2 analo}
M_t^n &:=& \log\frac{N^{\Lambda,n}(t \wedge\tau_{n_0}^n)}{v(t_n+t
\wedge
\tau_{n_0}^n)}\nonumber\\[-8pt]\\[-8pt]
&&{} + \int_0^{t\wedge\tau_{n_0}^n}
\biggl(\frac{\psi(N^{\Lambda,n}(r))}{N^{\Lambda,n}(r)} -
\frac{\psi(v(t_n+r))}{v(t_n+r)} +h^n(r) \biggr) \,dr,\nonumber
\end{eqnarray}
where $h^n$ is the drift compensator of $\log(N^{\Lambda,n})$
with respect to the filtration generated by the underlying $\Lambda
$-coalescent and where
\[
\tau_{n_0}^n:=\inf\{s> 0\dvtx N^{\Lambda,n}(s)\leq n_0\}.
\]
Then $M^n$ in (\ref{E M thm2 analo})
is a direct analogue of martingale (\ref{E M thm2}).
In particular, note that by definition of $t_n$, $M_0^n=0$, and as in
(\ref{E varM estim}),
\[
E((M_t^n)^2) \leq C t.
\]
Recall $\tau_{n_0}$ defined in (\ref{E_taun0}), and
note that with probability $1$, $\tau_{n_0}^n$ increases to $\tau_{n_0}$
as $n\to\infty$.
The arguments leading to Proposition \ref{Psummar}
apply in the current setting to yield for a fixed $\alpha^*<1/2$, and
for all 
$n$ (for $n\leq n_0$ the result holds trivially),
%
%
\begin{eqnarray}
\label{E exp est ncase}
P \biggl(\sup_{t\in[0,s]} \biggl|\log\frac{N^{\Lambda,n}(t\wedge
\tau
_{n_0}^n)}{v(t_n+t\wedge\tau_{n_0}^n)}
\biggr| \leq24 s^{\alpha^*} \biggr) &\geq&
1-O(s^{1-2 \alpha^*})\quad \mbox{and}
\\
\label{E exp est ncase moment}
P \biggl(\sup_{t\in[0,s]} \biggl|\log\frac{N^{\Lambda,n}(t\wedge
\tau
_{n_0}^n)}{v(t_n+t\wedge\tau_{n_0}^n)}
\biggr| \leq2(x+Cs) \biggr) &\geq&
1-O(r(x;s)).
\end{eqnarray}

Fix some subsequence $(n_k)_{k\geq1}$. 
We now show the first convergence statement in~(\ref{Esuffic}).
Choose any sequence $s_j$ of positive numbers decreasing to $0$ so that
%
%
\begin{equation}\label{E impo conv as}
\sum_j s_j^{1-2\alpha^*} < \infty.
\end{equation}
Next choose a further subsequence of $(n_k)_{k\geq1}$, denoted again
by $(n_j)_{j\geq1}$
to simplify notation, so that
%
%
\begin{eqnarray}
\label{E_vvN}
\lim_{j\to\infty} \int_0^{s_j} v(t_{n_j}+t) \,dt &=&
\infty,\nonumber\\[-8pt]\\[-8pt]
\lim_{j\to\infty} \frac{\int_{s_j}^s v(t_{n_j}+t) \,dt}{\int
_0^{s_j} v(t_{n_j}+t) \,dt} &=&
\lim_{j\to\infty} \frac{\int_{s_j}^s N^{\Lambda,n_j}(t) \,dt}{\int
_0^{s_j} N^{\Lambda,n_j}(t) \,dt} = 0,\nonumber
\end{eqnarray}
where the last limit is taken almost surely.
Note that
here we use observations (\ref{E v nonint}) and (\ref{E:abhl}) and the
following straightforward facts: for any fixed $0\leq a<b\leq s$,
$\int_a^b v(t_{n_j}+t) \,dt \uparrow\int_a^b v(t) \,dt$ and
$\int_a^b N^{\Lambda,n_j}(t) \,dt \uparrow\int_a^b N^{\Lambda}(t)
\,dt$. Due
to (\ref{E exp est ncase}), (\ref{E impo conv as}) and the
Borel--Cantelli lemma, we have
%
%
\begin{equation}
\label{Eas_conver}
\lim_{j\to\infty} \sup_{t\in[0,s_j]} \biggl|\frac{N^{\Lambda
,n_j}(t)}{v(t_{n_j}+t)}-1 \biggr|=0\qquad\mbox{almost surely}.
\end{equation}
The first statement in (\ref{Esuffic}) now follows by a simple
calculus fact:
if $(f_n)_{n\geq1}$, $(g_n)_{n\geq1}$, $f_n,g_n\dvtx[0,s]\to[0,\infty)$,
are two sequences of integrable functions such that for some positive sequence
$\delta_n\to0$ it is true that
\[
\lim_{n\to\infty} \int_0^{\delta_n} f_n(t) \,dt = \infty,\qquad
\lim_{n\to\infty} \frac{\int_{\delta_n}^s f_n(t) \,dt}{\int
_0^{\delta_n} f_n(t) \,dt} =
\lim_{n\to\infty} \frac{\int_{\delta_n}^s g_n(t) \,dt}{\int
_0^{\delta_n} g_n(t) \,dt} = 0,
\]
and
\[
\lim_{n\to\infty} \sup_{t\in[0,\delta_n]} \biggl|\frac
{f_n(t)}{g_n(t)}-1 \biggr|=0,
\]
then
\[
\lim_{n\to\infty} \frac{\int_0^s f_n(t) \,dt}{\int_0^s g_n(t) \,dt} =1.
\]

For the second convergence statement in (\ref{Esuffic}), note that
(similar to the argument
for Theorem \ref{Tsmalltime-mom}),
almost sure convergence (\ref{Eas_conver}) together with
estimate (\ref{E exp est ncase moment}) and the dominated convergence
theorem, yield
%
%
\begin{equation}
\label{E:unifrm convgce of E/v} \lim_{j\to\infty} \sup_{t\in
[0,s_j]} \biggl|\frac{E N^{\Lambda,n_j}(t)}{v(t_{n_j}+t)}-1
\biggr|=0\qquad
\mbox{almost surely}.
\end{equation}
Note that without loss of generality we may assume that
%
%
\begin{equation}
\label{E_ENEN}
\lim_{j\to\infty} \frac{\int_{s_j}^s EN^{\Lambda,n_j}(t)
\,dt}{\int
_0^{s_j} EN^{\Lambda,n_j}(t) \,dt} = 0.
\end{equation}
The previous argument applies.

The final statement of Theorem \ref{C:extraasym} will follow from
Corollary \ref{C:extraasym as},
which is stated and proved in next subsection.

\subsubsection{Discussion on almost sure convergence}
It is an open question whether the convergence of Theorem \ref
{C:extraasym} holds almost surely.
Our technique seems too crude to verify it in general, yet we
offer below a partial result in this direction.
One standard approach would be to use the monotonicity
\[
\int_0^s N^{\Lambda,n}(t) \,dt
\leq\int_0^s N^{\Lambda,n+1}(t) \,dt \quad\mbox{and}\quad
\int_0^s v(t_{n}+t) \,dt
\leq\int_0^s v(t_{n+1}+t) \,dt.
\]
It would suffice to find a subsequence $n_j$ along which convergence holds
in the almost sure sense, and in addition, such that
%
%
\begin{equation}
\label{E extra for as}
\lim_{j\to\infty}
\frac{\int_0^s v(t_{n_j}+t) \,dt}{\int_0^s v(t_{n_{j+1}}+t) \,dt}=1.
\end{equation}
\begin{coro}
\label{C:extraasym as}
Assume that $\alpha^*<1/2$ is fixed, and that
two sequences $(s_j)_{j\geq1}$ and $(n_j)_{j\geq1}$ are given where
$n_j$ is nondecreasing.
If in addition to (\ref{E impo conv as})
and (\ref{E extra for as}), we have
%
%
\begin{eqnarray}
\label{E_vv}
\lim_{j\to\infty} \int_0^{s_j} v(t_{n_j}+t) \,dt &=& \infty,\qquad
\lim_{j\to\infty} \frac{\int_{s_j}^s v(t_{n_j}+t) \,dt}{\int
_0^{s_j} v(t_{n_j}+t) \,dt} = 0\quad
\mbox{and}
\\
\label{E vv a}
\lim_{j\to\infty} \frac{\int_{s_j}^s v(t) \,dt}{\int_{s_j}^s
v(t_{n_j}+t) \,dt}&<&\infty,
\end{eqnarray}
then the convergence of Theorem \ref{C:extraasym} holds almost surely.
\end{coro}
\begin{pf}
As discussed above, due to (\ref{E extra for as}) and monotonicity,
it suffices to show convergence as stated in Theorem
\ref{C:extraasym} along the sequence $(n_j)_{j\geq1}$. Due to the
Borel--Cantelli lemma, (\ref{E exp est ncase}), (\ref{E impo conv
as}), (\ref{E:unifrm convgce of E/v}) and the fact
\[
P\Bigl(\limsup_j \{\tau_{n_0}^{n_j}< s_j\}\Bigr)=0,
\]
we have, as for Theorems \ref{Tsmalltime} and \ref{Tsmalltime-mom}, that
\[
\lim_{j\to\infty} \frac{\int_0^{s_j} N^{\Lambda,n_j}(t)
\,dt}{\int
_0^{s_j} v(t_{n_j}+t) \,dt}=1\qquad\mbox{almost surely}
\]
and
\[
\lim_{j\to\infty} \frac{\int_0^{s_j} EN^{\Lambda,n_j}(t)
\,dt}{\int
_0^{s_j} v(t_{n_j}+t) \,dt}=1.
\]
Due to (\ref{E_vv}),
we have
\[
\liminf_{j\to\infty} \frac{\int_0^s N^{\Lambda,n_j}(t) \,dt}{\int
_0^s v(t_{n_j}+t) \,dt}\geq1
\qquad\mbox{almost surely}
\]
and
\[
\liminf_{j\to\infty} \frac{\int_0^s EN^{\Lambda,n_j}(t)
\,dt}{\int
_0^s v(t_{n_j}+t) \,dt}\geq1.
\]
For the corresponding upper bound on the $\limsup$, note that due to Theorem
\ref{Tsmalltime} (resp. Theorem \ref{Tsmalltime-mom}) there exists a
positive finite random variable
$C_0$ (resp. positive constant $C_0$) such that
\[
\frac{\int_{s_j}^s N^{\Lambda}(t) \,dt}{\int_{s_j}^s v(t) \,dt}
\leq
1+C_0\qquad\mbox{a.s.},
\biggl(\mbox{resp. }\frac{\int_{s_j}^s EN^{\Lambda}(t) \,dt}{\int
_{s_j}^s v(t) \,dt} \leq1+C_0 \biggr) ,
\mbox{ for all }j\geq1.
\]
Due to (\ref{E_vv}) and (\ref{E vv a}) and monotonicity $N^{\Lambda
,n_j}(t)\leq N^{\Lambda}(t)$ (with probability $1$), we now have both
\[
\lim_j \frac{\int_{s_j}^s N^{\Lambda,n_j}(t) \,dt}{\int
_0^{s_j}v(t_{n_j}+t) \,dt} = 0\qquad
\mbox{almost surely},\quad\mbox{and}\quad
\lim_j \frac{\int_{s_j}^s EN^{\Lambda,n_j}(t) \,dt}{\int
_0^{s_j}v(t_{n_j}+t) \,dt} = 0,
\]
which completes the argument.
\end{pf}

Taking for example
$\alpha^*=1/4$, $s_j=1/j^3$,
and $n_j=\exp(\log^2{j})$ (resp. $n_j=j^{\eta}$  with $\eta>3 (\alpha-1)$)
in the case of Kingman (resp. Beta) coalescent, one can verify (left to
the reader)
the hypotheses of the last corollary, implying the final statement of
Theorem \ref{C:extraasym}.

\begin{appendix}\label{ap:binom}
\section*{Appendix: Binomial calculations}

\begin{lemma}
\label{Lbincalc} If $X$ has $\operatorname{Binomial}(n,p)$ distribution and if $Y=X -
\mathbf{1}_{\{X>0\}}$, then:
%
%
\begin{eqnarray}\label{Esecmom}
\mbox{\textup{(i) }\hspace*{15pt}} E Y &=& np - 1 + (1-p)^n;\nonumber\\
\mbox{\textup{(ii) }} \operatorname{var} (Y) &=&
np (1-p) + (1-p)^n \bigl(1-(1-p)^n\bigr) - 2np (1-p)^n;\\
\mbox{\textup{(iii)}\hspace*{13pt}} E Y^2 &=& -np - n p^2 + n^2 p^2 + 1 - (1-p)^n.\nonumber
\end{eqnarray}
\end{lemma}
\begin{pf}
Property (i) is trivial, (ii) follows easily from the fact that
\[
\operatorname{cov}\bigl(X, \mathbf{1}_{\{X>0\}}\bigr) = np (1-p)^n
\]
and (iii) is implied by (i) and (ii).
\end{pf}
\begin{coro}
\label{C:secmombd} If $X$ has $\operatorname{Binomial}(n,1-p)$ distribution and if $Y=
X + \mathbf{1}_{\{X<n\}}$, then
\[
E \biggl[ \biggl(\frac{n-Y}{n} \biggr)^2 \biggr] = O(p^2).
\]
\end{coro}
\begin{pf}
Note that $n-Y$ has the distribution of the variable $Y$ from
Lemma~\ref{Lbincalc}. Hence its second moment is given in (\ref{Esecmom}).
Since for $p<1/n$ we have
\[
(1-p)^n = 1-np + O(n^2p^2),
\]
the claim of the corollary is true in this case. Now if $p\geq1/n$
then $np = O(n^2 p^2)$ therefore the largest term in (\ref{Esecmom})
is again of order $n^2 p^2$.
\end{pf}
\begin{lemma}
\label{Lbinlogcalc} There exists $n_0\in\mathbb{N}$ and $C_0<\infty
$ such
that for all $n\geq n_0$ and all $p\leq1/4$, if $X$ has
$\operatorname{Binomial}(n,1-p)$ distribution, then
%
\[
\biggl| E\bigl[ \log\bigl(X+\mathbf{1}_{\{X<n\}}\bigr) - \log{n}\bigr] + \frac{np-1+(1-p)^n}{n}
\biggr| \leq C_0 p^2
\]
and
\[
E\bigl[\bigl( \log\bigl(X+\mathbf{1}_{\{X<n\}}\bigr) - \log{n}\bigr)^2\bigr] \leq C_0 p^2.
\]
\end{lemma}
\begin{pf}
Let $Y=n-X$ as before, and abbreviate
%
%
\begin{equation}
\label{ET}
T\equiv T_n := \log\bigl(X+\mathbf{1}_{\{X<n\}}\bigr) - \log{n}=
\log\biggl(1-
\frac{Y-\mathbf{1}_{\{Y>0\}}}{n} \biggr).
\end{equation}
We split the computation according to the event
\[
A_n =\{Y \leq n/2\},
\]
whose complement due to a large deviation bound has probability
bounded by
%
%
\begin{equation}
\label{ELDPbound} \exp\biggl\{ -n
\biggl(\frac{1}{2}\log{\frac{1}{2p}}+
\frac{1}{2}\log{\frac{1}{2(1-p)}} \biggr) \biggr\} = 2^n p^{n/2}
(1-p)^{n/2},
\end{equation}
uniformly in $p\leq1/4$ and $n$.
On $A_n^c$ we have $|T| \leq\log{n}$, and on $A_n$ we apply a
calculus fact,
%
%
\begin{equation}\label{Ecalc_fact}
|{\log}(1-x) + x| \leq\frac{x^2}{2(1-x)}\leq x^2,\qquad x\in[0,1/2],
\end{equation}
to obtain
%
\[
\biggl|E[T] + E \biggl[\frac{Y-\mathbf{1}_{\{Y>0\}}}{n} \mathbf
{1}_{A_n} \biggr]
\biggr| \leq(\log{n}) P(A_n^c) +
E \biggl[\frac{(Y-\mathbf{1}_{\{Y>0\}})^2}{n^2} \mathbf
{1}_{A_n} \biggr].
\]
Furthermore, since $(Y-\mathbf{1}_{\{Y>0\}})/n\leq1$, we conclude
%
%
\begin{equation}
\label{Ehelping}\qquad
\biggl|E[T] + E \biggl[\frac{Y-\mathbf{1}_{\{Y>0\}}}{n} \biggr] \biggr|
\leq
(\log{n}+1) P(A_n^c) + E \biggl[\frac{(Y-\mathbf{1}_{\{Y>0\}
})^2}{n^2} \biggr].
\end{equation}
Note that by Corollary \ref{C:secmombd} and Lemma \ref{Lbincalc}(i),
in order to prove the first estimate of the lemma, it remains to
show
%
%
\begin{equation}
\label{EPlognbound}
(\log{n}) P(A_n^c)\leq(\log{n}) 2^n p^{n/2}
(1-p)^{n/2}\leq C p^2
\end{equation}
for some $C<\infty$, all $p\in[0,1/4]$, and all $n$ large. Now
consider $f\dvtx p\mapsto(p(1-p))^{n/2}/p^2$. Its derivative at $p$
equals $g(p) (n(1-2p)/2 - 2(1-p))$ where $g(p)$ is a positive
function. It is easy to check that if $p\leq1/4$, then $n(1-2p)/2 -
2p^2(1-p)>0$ for all $n\geq6$. Therefore $f$ is an increasing
function of $p$,
so in order to verify (\ref{EPlognbound})
for all $p\leq1/4$,
it suffices to check it for $p=1/4$.
This corresponds to having
$(\log{n}) 2^n (3/16)^{n/2} \leq C/16 $,
that will hold for all large $n=n(C)$ given a $C > 0$.

For the second estimate, again use the partitioning according to
$A_n$ and (\ref{Ecalc_fact}) to obtain
\[
ET^2 \leq E[T^2 \mathbf{1}_{A_n}]+ \log^2{n} P(A_n^c) \leq
\biggl(\frac{3}{2} \biggr)^2
E \biggl[\frac{(Y-\mathbf{1}_{\{Y>0\}})^2}{n^2} \biggr] + \log^2{n}
P(A_n^c),
\]
which differs from (\ref{Ehelping}) only by an extra factor of order
$\log{n}$ multiplying
$P(A_n^c)$, so the previous argument carries over.
\end{pf}
%
%
\begin{lemma}
\label{Lbinlogcalc1} There exists $n_0\in\mathbb{N}$ and $K_0<\infty
$ such
that for all $n\geq n_0$, $p\leq1/4$ and $c>0$, if $X$ has
$\operatorname{Binomial}(n,1-p)$ distribution, then
\[
E\bigl[ \exp\bigl\{c\bigl[\log\bigl(X+\mathbf{1}_{\{X<n\}}\bigr) - \log{n}\bigr]^2\bigr\} -1\bigr] \leq e^{9c/4}
K_0 p^2.
\]
\end{lemma}
\begin{pf}
The strategy is the same as that used for the second estimate in the
previous lemma, some details
are left to the reader.

Recall that $Y=n-X$ and observe that
\begin{eqnarray*}
E[e^{c T^2}-1] &\leq& n^{c \log{n}} P(A_n^c)
+ E[(e^{c T^2}-1)\mathbf{1}_{A_n}]\\
&\leq&
n^{c \log{n}} P(A_n^c) +
E \biggl[ \biggl(\exp\biggl\{c \frac{9(Y-\mathbf{1}_{\{Y>0\}
})^2}{4n^2} \biggr\}-1 \biggr)\mathbf{1}_{A_n} \biggr]\\
&\leq&
n^{c \log{n}} P(A_n^c) +
E \biggl[\exp\biggl\{c \frac{9(Y-\mathbf{1}_{\{Y>0\}})^2}{4n^2}
\biggr\}
-1 \biggr].
\end{eqnarray*}
Hence it suffices to show that for some $K_0$, all $c>0$ and all
$n$, $p$ as specified above, we have
%
%
\begin{equation}\label{Esufi}
E \biggl[\exp\biggl\{c \frac{(Y-\mathbf{1}_{\{Y>0\}})^2}{n^2} \biggr\}
-1 \biggr]
\leq e^c K_0 p^2.
\end{equation}
Without loss of generality, one can assume that $c>1$.

The left-hand side above
\[
\sum_{k=1}^n \pmatrix{n \cr k} p^k (1-p)^{n-k}
\bigl(e^{c(k-1)^2/n^2}-1 \bigr)
\]
can be bounded, using  Taylor's expansion, by
%
%
\begin{eqnarray}\label{Esufi1}\quad
&&\sum_{k=1}^n \pmatrix{n \cr k} p^k (1-p)^{n-k}
\biggl\{ c\frac{(k-1)^2}{n^2}
+ \frac{e^c}{2} \frac{(k-1)^4}{n^4} \biggr\}\nonumber\\
&&\qquad=\frac{c}{n^2} \bigl(E(Y-1)^2 - P(Y=0)\bigr)\\
&&\qquad\quad{} + \frac{e^c}{2n^4} \bigl(E(Y-1)^4 -
P(Y=0)\bigr).\nonumber
\end{eqnarray}
Next compute
%
%
\begin{eqnarray}\label{Elast}
&&E(Y-1)^2 - P(Y=0)\nonumber\\
&&\qquad=
\operatorname{Var}(Y-1) +\bigl( E(Y-1)\bigr)^2- P(Y=0) \nonumber\\
&&\qquad=
np(1-p)+(np-1)^2 -(1-p)^n\qquad\mbox{[recall (\ref
{Eexppowbd})]}\\
&&\qquad\leq
np(1-p)+(np-1)^2 - e^{-np} + 2np^2/3
\nonumber\\
&&\qquad\leq
(np)^2 + O(np^2),\nonumber
\end{eqnarray}
where, for the last inequality, we recall that $e^{-x}-1 +x>0$ for $x>0$.
Similarly, using the fact
\[
(y-1)^4= y(y-1)(y-2)(y-3)+2y(y-1)(y-2)+(y-1)^2
\]
as well as the expressions for Binomial factorial moments, we have
%
%
\begin{eqnarray}
E(Y &-& 1)^4 - P(Y=0)\\
&=&
n(n-1)(n-2)(n-3)p^4\nonumber\\
&&{}+2n(n-1)(n-2)p^3+E(Y-1)^2-P(Y=0)
\nonumber\\
&\leq&
n^4 p^4 + 2 n^3 p^3 + (np)^2 + O(np^2)\nonumber\\
\label{Elasta}
&\leq&
4 n^4 p^2 + O(np^2).
\end{eqnarray}
Now (\ref{Esufi1})--(\ref{Elasta}) yield (\ref{Esufi}), and therefore
the statement of the lemma, with
appropriately chosen $n_0$.
\end{pf}
\end{appendix}

\section*{Acknowledgments}
We are grateful to Etienne Pardoux
for his careful reading of a preliminary draft,
and to anonymous referees for their valuable comments and suggestions.
The research project on which both
the current work and \cite{bbl2} report,
started at the time when
J. Berestycki was a junior faculty at Universit\'{e} de Provence,
while
N. Berestycki was a postdoctoral fellow and V. Limic a junior faculty
at the University of British Columbia.

%

%
\printaddresses


\begin{thebibliography}{29}

\bibitem{abhl}
%
\begin{bmisc}[unstr]
\bauthor{\bsnm{Angel},~\bfnm{O.}\binits{O.}},
\bauthor{\bsnm{Berestycki},~\bfnm{N.}\binits{N.}},
\bauthor{\bsnm{Hammond},~\bfnm{A.}\binits{A.}} \AND
\bauthor{\bsnm{Limic},~\bfnm{V.}\binits{V.}}
(\byear{2009}).
\bhowpublished{Global divergence of spatial coalescents. In preparation.}
\end{bmisc}
%
\endbibitem

\bibitem{barlowetal}
%
\begin{barticle}[mr]
\bauthor{\bsnm{Barlow},~\bfnm{M.~T.}\binits{M.~T.}},
\bauthor{\bsnm{Jacka},~\bfnm{S.~D.}\binits{S.~D.}} \AND
\bauthor{\bsnm{Yor},~\bfnm{M.}\binits{M.}}
(\byear{1986}).
\btitle{Inequalities for a pair of processes stopped at a random time}.
\bjournal{Proc. London Math. Soc.}
\bvolume{52}
\bpages{142--172}.
\bid{mr={812449}}
\end{barticle}
%
\endbibitem

\bibitem{bbl2}
%
\begin{bmisc}[unstr]
\bauthor{\bsnm{Berestycki},~\bfnm{J.}\binits{J.}},
\bauthor{\bsnm{Berestycki},~\bfnm{N.}\binits{N.}} \AND
\bauthor{\bsnm{Limic},~\bfnm{V.}\binits{V.}}
(\byear{2008}).
\bhowpublished{Interpreting
$\Lambda$-coalescent speed of coming down from infinity via particle
representation of super-processes. In preparation.}
\end{bmisc}
%
\endbibitem

\bibitem{bbs2}
%
\begin{barticle}[mr]
\bauthor{\bsnm{Berestycki},~\bfnm{Julien}\binits{J.}},
\bauthor{\bsnm{Berestycki},~\bfnm{Nathana{\"e}l}\binits{N.}} \AND
\bauthor{\bsnm{Schweinsberg},~\bfnm{Jason}\binits{J.}}
(\byear{2007}).
\btitle{Beta-coalescents and continuous stable random trees}.
\bjournal{Ann. Probab.}
\bvolume{35}
\bpages{1835--1887}.
\bid{mr={2349577}}
\end{barticle}
%
\endbibitem

\bibitem{bbs1}
%
\begin{barticle}[mr]
\bauthor{\bsnm{Berestycki},~\bfnm{Julien}\binits{J.}},
\bauthor{\bsnm{Berestycki},~\bfnm{Nathana{\"e}l}\binits{N.}} \AND
\bauthor{\bsnm{Schweinsberg},~\bfnm{Jason}\binits{J.}}
(\byear{2008}).
\btitle{Small-time behavior of beta coalescents}.
\bjournal{Ann. Inst. H. Poincar\'e Probab. Statist.}
\bvolume{44}
\bpages{214--238}.
\bid{mr={2446321}}
\end{barticle}
%
\endbibitem

\bibitem{blg3}
%
\begin{barticle}[mr]
\bauthor{\bsnm{Bertoin},~\bfnm{Jean}\binits{J.}} \AND
\bauthor{\bsnm{Le~Gall},~\bfnm{Jean-Francois}\binits{J.-F.}}
(\byear{2006}).
\btitle{Stochastic flows associated to coalescent processes. {III}. {L}imit
theorems}.
\bjournal{Illinois J. Math.}
\bvolume{50}
\bpages{147--181}.
\bid{mr={2247827}}
\end{barticle}
%
\endbibitem

\bibitem{bosz98}
%
\begin{barticle}[mr]
\bauthor{\bsnm{Bolthausen},~\bfnm{E.}\binits{E.}} \AND
\bauthor{\bsnm{Sznitman},~\bfnm{A.-S.}\binits{A.-S.}}
(\byear{1998}).
\btitle{On {R}uelle's probability cascades and an abstract cavity method}.
\bjournal{Comm. Math. Phys.}
\bvolume{197}
\bpages{247--276}.
\bid{mr={1652734}}
\end{barticle}
%
\endbibitem

\bibitem{darnor}
%
\begin{barticle}[mr]
\bauthor{\bsnm{Darling},~\bfnm{R.~W.~R.}\binits{R.~W.~R.}} \AND
\bauthor{\bsnm{Norris},~\bfnm{J.~R.}\binits{J.~R.}}
(\byear{2008}).
\btitle{Differential equation approximations for {M}arkov chains}.
\bjournal{Probab. Surv.}
\bvolume{5}
\bpages{37--79}.
\bid{mr={2395153}}
\end{barticle}
%
\endbibitem

\bibitem{ddsj}
%
\begin{barticle}[mr]
\bauthor{\bsnm{Delmas},~\bfnm{Jean-Fran{\c{c}}ois}\binits{J.-F.}},
\bauthor{\bsnm{Dhersin},~\bfnm{Jean-St{\'e}phane}\binits{J.-S.}} \AND
\bauthor{\bsnm{Siri-Jegousse},~\bfnm{Arno}\binits{A.}}
(\byear{2008}).
\btitle{Asymptotic results on the length of coalescent trees}.
\bjournal{Ann. Appl. Probab.}
\bvolume{18}
\bpages{997--1025}.
\bid{mr={2418236}}
\end{barticle}
%
\endbibitem

\bibitem{dong}
%
\begin{barticle}[unstr]
\bauthor{\bsnm{Dong},~\bfnm{Rui}\binits{R.}},
\bauthor{\bsnm{Gnedin},~\bfnm{Alexander}\binits{A.}} \AND
\bauthor{\bsnm{Pitman},~\bfnm{Jim}\binits{J.}}
(\byear{2007}).
\btitle{Exchangeable partitions derived from Markovian coalescents}.
\bjournal{Ann. Appl. Probab.}
\bvolume{17}
\bpages{1172--1201}.
\bid{mr={2344303}}
\end{barticle}
%
\endbibitem

\bibitem{durrett}
%
\begin{bbook}[mr]
\bauthor{\bsnm{Durrett},~\bfnm{Richard}\binits{R.}}
(\byear{2004}).
\btitle{Probability: Theory and Examples}, \bedition{3rd} ed.
\bpublisher{Duxbury Press}, \baddress{Belmont, CA}.
\bid{mr={1609153}}
\end{bbook}
%
\endbibitem

\bibitem{durrettDNA}
%
\begin{bbook}[mr]
\bauthor{\bsnm{Durrett},~\bfnm{Rick}\binits{R.}}
(\byear{2002}).
\btitle{Probability Models for {DNA} Sequence Evolution}.
\bpublisher{Springer}, \baddress{New York}.
\bid{mr={1903526}}
\end{bbook}
%
\endbibitem

\bibitem{eldwak}
%
\begin{barticle}[unstr]
\bauthor{\bsnm{Eldon},~\bfnm{B.}\binits{B.}} \AND
\bauthor{\bsnm{Wakeley},~\bfnm{J.}\binits{J.}}
(\byear{2006}).
\btitle{Coalescent processes when the
distribution of offspring number among individuals is highly skewed}.
\bjournal{Genetics}
\bvolume{172}
\bpages{2621--2633}.
\end{barticle}
%
\endbibitem

\bibitem{ewens}
%
\begin{bbook}[mr]
\bauthor{\bsnm{Ewens},~\bfnm{Warren~J.}\binits{W.~J.}}
(\byear{2004}).
\btitle{Mathematical Population Genetics}.
\bpublisher{Springer}, \baddress{New York}.
\bid{mr={554616}}
\end{bbook}
%
\endbibitem

\bibitem{hoeffding}
%
\begin{barticle}[mr]
\bauthor{\bsnm{Hoeffding},~\bfnm{Wassily}\binits{W.}}
(\byear{1963}).
\btitle{Probability inequalities for sums of bounded random variables}.
\bjournal{J.~Amer. Statist. Assoc.}
\bvolume{58}
\bpages{13--30}.
\bid{mr={0144363}}
\end{barticle}
%
\endbibitem

\bibitem{king82}
%
\begin{barticle}[mr]
\bauthor{\bsnm{Kingman},~\bfnm{J.~F.~C.}\binits{J.~F.~C.}}
(\byear{1982}).
\btitle{The coalescent}.
\bjournal{Stochastic Process. Appl.}
\bvolume{13}
\bpages{235--248}.
\bid{mr={671034}}
\end{barticle}
%
\endbibitem

\bibitem{king82b}
%
\begin{barticle}[mr]
\bauthor{\bsnm{Kingman},~\bfnm{J.~F.~C.}\binits{J.~F.~C.}}
(\byear{1982}).
\btitle{On the genealogy of large populations}.
\bjournal{J. Appl. Probab.}
\bvolume{19A}
\bpages{27--43}.
\bid{mr={633178}}
\end{barticle}
%
\endbibitem

\bibitem{ls}
%
\begin{barticle}[mr]
\bauthor{\bsnm{Limic},~\bfnm{Vlada}\binits{V.}} \AND
\bauthor{\bsnm{Sturm},~\bfnm{Anja}\binits{A.}}
(\byear{2006}).
\btitle{The spatial {$\Lambda$}-coalescent}.
\bjournal{Electron. J. Probab.}
\bvolume{11}
\bpages{363--393}.
\bid{mr={2223040}}
\end{barticle}
%
\endbibitem

\bibitem{marine}
%
\begin{barticle}[unstr]
\bauthor{\bsnm{Li},~\bfnm{G.}\binits{G.}} \AND
\bauthor{\bsnm{Hedgecock},~\bfnm{D.}\binits{D.}}
(\byear{1998}).
\btitle{Genetic heterogeneity, detected
by PCR SSCP, among samples of larval Pacific oysters (Crassostrea gigas)
supports the hypothesis of large variance in reproductive success}.
\bjournal{Canad. J. Fish. Aquat. Sci.}
\bvolume{55}
\bpages{1025--1033}.
\end{barticle}
%
\endbibitem

\bibitem{mohle}
%
\begin{barticle}[mr]
\bauthor{\bsnm{M{\"o}hle},~\bfnm{M.}\binits{M.}}
(\byear{2006}).
\btitle{On sampling distributions for coalescent processes with simultaneous
multiple collisions}.
\bjournal{Bernoulli}
\bvolume{12}
\bpages{35--53}.
\bid{mr={2202319}}
\end{barticle}
%
\endbibitem

\bibitem{morrisetal}
%
\begin{barticle}[unstr]
\bauthor{\bsnm{Morris},~\bfnm{A. P.}\binits{A. P.}},
\bauthor{\bsnm{Whittaker},~\bfnm{J. C.}\binits{J. C.}} \AND
\bauthor{\bsnm{Balding},~\bfnm{D. J.}\binits{D. J.}}
(\byear{2002}).
\btitle{Fine-scale mapping of disease loci via
shattered coalescent modeling of genealogies}.
\bjournal{Amer. J. Hum. Genet.}
\bvolume{70}
\bpages{686--707}.
\end{barticle}
%
\endbibitem

\bibitem{pit99}
%
\begin{barticle}[mr]
\bauthor{\bsnm{Pitman},~\bfnm{Jim}\binits{J.}}
(\byear{1999}).
\btitle{Coalescents with multiple collisions}.
\bjournal{Ann. Probab.}
\bvolume{27}
\bpages{1870--1902}.
\bid{mr={1742892}}
\end{barticle}
%
\endbibitem

\bibitem{revuzyor}
%
\begin{bbook}[mr]
\bauthor{\bsnm{Revuz},~\bfnm{Daniel}\binits{D.}} \AND
\bauthor{\bsnm{Yor},~\bfnm{Marc}\binits{M.}}
(\byear{1994}).
\btitle{Continuous Martingales and {B}rownian Motion},
\bedition{2nd} ed.
\bseries{Grundlehren der Mathematischen Wissenschaften [Fundamental Principles
of Mathematical Sciences]}
\bvolume{293}.
\bpublisher{Springer}, \baddress{Berlin}.
\bid{mr={1303781}}
\end{bbook}
%
\endbibitem

\bibitem{rogerswilliams}
%
\begin{bbook}[mr]
\bauthor{\bsnm{Rogers},~\bfnm{L.~C.~G.}\binits{L.~C.~G.}} \AND
\bauthor{\bsnm{Williams},~\bfnm{David}\binits{D.}}
(\byear{1987}).
\btitle{Diffusions, Markov Processes and Martingales},
\bedition{2nd} ed.
\bpublisher{Cambridge Univ. Press}, \baddress{Cambridge}.
\bid{mr={1780932}}
\end{bbook}
%
\endbibitem

\bibitem{delapena}
%
\begin{barticle}[mr]
\bauthor{\bparticle{de~la }\bsnm{Pe{\~n}a},~\bfnm{Victor~H.}\binits{V.~H.}}
(\byear{1999}).
\btitle{A general class of exponential inequalities for martingales and
ratios}.
\bjournal{Ann. Probab.}
\bvolume{27}
\bpages{537--564}.
\bid{mr={1681153}}
\end{barticle}
%
\endbibitem

\bibitem{sag99}
%
\begin{barticle}[mr]
\bauthor{\bsnm{Sagitov},~\bfnm{Serik}\binits{S.}}
(\byear{1999}).
\btitle{The general coalescent with asynchronous mergers of ancestral lines}.
\bjournal{J.~Appl. Probab.}
\bvolume{36}
\bpages{1116--1125}.
\bid{mr={1742154}}
\end{barticle}
%
\endbibitem

\bibitem{sj}
%
\begin{bmisc}[unstr]
\bauthor{\bsnm{Siri-Jegousse},~\bfnm{A.}\binits{A.}}
(\byear{2008}).
\bhowpublished{\textit{Etude de Processus Appliqu\'{e}s
Aux Mod\'{e}les de Population}. Ph.D. thesis, Paris V.}
\end{bmisc}
%
\endbibitem

\bibitem{sch1}
%
\begin{barticle}[mr]
\bauthor{\bsnm{Schweinsberg},~\bfnm{Jason}\binits{J.}}
(\byear{2000}).
\btitle{A necessary and sufficient condition for the {$\Lambda
$}-coalescent to
come down from infinity}.
\bjournal{Electron. Comm. Probab.}
\bvolume{5}
\bpages{1--11}.
\bid{mr={1736720}}
\end{barticle}
%
\endbibitem\

\bibitem{sch2}
%
\begin{barticle}[mr]
\bauthor{\bsnm{Schweinsberg},~\bfnm{Jason}\binits{J.}}
(\byear{2000}).
\btitle{Coalescents with simultaneous multiple collisions}.
\bjournal{Electron. J.~Probab.}
\bvolume{5}
\bpages{50}.
\bid{mr={1781024}}
\end{barticle}
%
\endbibitem

\end{thebibliography}
\end{document}